\documentclass[preprint,10pt]{elsarticle}
\usepackage{graphicx}
\usepackage[ruled,vlined]{algorithm2e}
\usepackage{algorithmic}
\usepackage{epsfig}
\usepackage{epstopdf}
\usepackage{amssymb,amsmath}
\usepackage{amsmath}
\usepackage{amsthm}
\usepackage{amsfonts}
\usepackage{psfrag}
\usepackage{rotating}
\usepackage{latexsym}
\usepackage{dsfont}
\usepackage{stmaryrd}
\usepackage{amsthm}
\usepackage{amssymb}
\usepackage{amsmath}
\usepackage{mathtools}
\usepackage{mathrsfs}   
\usepackage{siunitx}
\usepackage{pgfpages}
\usepackage[Euler]{upgreek} 
\usepackage{isomath}
\usepackage{physics}
\usepackage{pifont} 
\usepackage{bbm}
\usepackage{empheq}
\usepackage[thicklines]{cancel}
\usepackage{subcaption}
\usepackage{natbib}

\usepackage[force,almostfull]{textcomp}
\usepackage{fancyvrb}
\usepackage{textcomp}
\usepackage{url}
\usepackage{ifdraft}
\usepackage{url}
\usepackage{multirow}
\usepackage[colorlinks=true,linkcolor=blue,citecolor=green,urlcolor=red]{hyperref}
\usepackage{rotating}
\usepackage{xspace}
\usepackage{array}
\usepackage{multido}
\usepackage{enumerate}
\usepackage[utf8]{inputenc}
\newlength\figureheight 
\newlength\figurewidth 
\usepackage{graphics}
\usepackage{pgfplots}
\pgfplotsset{compat=newest}
\pgfplotsset{plot coordinates/math parser=false}
\allowdisplaybreaks

\usepackage{scalefnt}

\newtheoremstyle{specialcasestyle}{1mm}{1mm}{\upshape}{}{\bfseries\upshape}{.}{0mm}{}
\theoremstyle{specialcasestyle}

\newtheorem{assump}{Assumption}

\newtheorem{rem}{Remark}
\newtheorem{thm}{Theorem}

\theoremstyle{remark}          

\newcommand*{\Var}[1]{\ensuremath{\mathrm{Var}\left[#1\right]}}
\newcommand*{\E}[1]{\ensuremath{\mathbb{E}\left[#1\right]}}
\newcommand*{\prob}[1]{\ensuremath{\mathbb{P}\left[#1\right]}}
\newcommand*{\tol}{\ensuremath{\mathrm{TOL}}}

\usepackage{amssymb}
\usepackage{amsmath}

\journal{Journal of Computational and Applied Mathematics}

\begin{document}

\begin{frontmatter}



\title{Multi-index importance sampling for McKean--Vlasov stochastic differential equations}


\author[label2]{Nadhir Ben Rached}
\affiliation[label2]{organization={Department of Statistics},
            addressline={School of Mathematics, University of Leeds}, 
            city={Leeds},
            postcode={LS2 9JT}, 
            country={United Kingdom}}
            
\author[label3]{Abdul-Lateef Haji-Ali}
\affiliation[label3]{organization={Department of Actuarial Mathematics and Statistics},
            addressline={School of Mathematical and Computer Sciences, Heriot-Watt University}, 
            city={Edinburgh},
            postcode={EH14 4AS}, 
            country={United Kingdom}}
            
\author[label1]{Shyam Mohan Subbiah Pillai \corref{cor1}} 
\cortext[cor1]{Corresponding author, \texttt{shyam.pillai@kaust.edu.sa}}
\affiliation[label1]{organization={Computer, Electrical and Mathematical Sciences \& Engineering Division (CEMSE)},
            addressline={King Abdullah University of Science and Technology (KAUST)}, 
            city={Thuwal},
            country={Saudi Arabia}}

\author[label1]{Ra\'ul Tempone}            

\begin{abstract}
This work addresses the estimation of rare-event quantities expressed as expectations of smooth observables of solutions to a broad class of McKean--Vlasov stochastic differential equations (MV-SDEs). Building on the double loop Monte Carlo (DLMC) method with stochastic optimal control-based importance sampling (IS) introduced by \cite{my_dlmcis}, this work extends this framework to the multi-index Monte Carlo (MIMC) setting. The resulting multi-index DLMC estimator mitigates the explosion of the coefficient of variation for rare event quantities. Moreover, it exploits the sampling efficiency of MIMC by leveraging the propagation of chaos to ensure mixed-difference variances vanish in the mean-field limit. The complexity analysis relies on assumptions on mixed-difference bias and variance decay, similar to standard MIMC assumptions. Although not rigorously proved, this work presents strong numerical evidence in support of these assumptions. The primary contribution of this work is the novel numerical integration of the MIMC method with IS for MV-SDEs. This approach reduces the computational complexity from $\order{\tol_{\mathrm{r}}^{-4}}$ for the DLMC estimator to $\order{\tol_{\mathrm{r}}^{-2} (\log \tol_{\mathrm{r}}^{-1})^2}$, enabling an accurate estimation of rare-event quantities within a prescribed relative error tolerance $\tol_{\mathrm{r}}$. Numerical experiments on the Kuramoto model from statistical physics demonstrate computational savings of several orders of magnitude for the multi-index DLMC estimator with IS, compared with the standard Monte Carlo (MC) method.
\end{abstract}

%

\begin{keyword}
McKean--Vlasov stochastic differential equation \sep importance sampling \sep multi-index Monte Carlo \sep decoupling approach \sep double loop Monte Carlo.

\MSC[2020] 60H35 \sep 65C05 \sep 65C30 \sep 65C35.

\end{keyword}

\end{frontmatter}


\section{Introduction}
\label{sec:intro}

McKean--Vlasov stochastic differential equations (MV-SDEs) are a class of SDEs whose drift and diffusion depend on the law of the solution~\citep{mckean_vlasov}. These SDEs arise as mean-field limits of stochastic interacting particle systems with applications in  pedestrian modeling~\citep{pedestrian_app}, animal behavior~\citep{animal_app}, biology~\citep{biology_app}, finance~\citep{finance_app}, and chemistry~\citep{chemical_app}. This work focuses on building a computationally efficient Monte Carlo (MC) estimator for $\E{G(X(T))}$, where $X:[0, T] \cross \Omega \rightarrow \mathbb{R}^d$ denotes the McKean--Vlasov process, and $G:\mathbb{R}^d \rightarrow \mathbb{R}$ represents a sufficiently smooth rare-event observable, meaning $\E{G(X(T))}$ is dominated by contributions from the tail of the law of $X$. For example, $G(x)$ can be a mollified version of $\mathbbm{1}_{\{x>K\}}$, for some large threshold $K$. The standard MC method is prohibitively inefficient in this regime because the probability of sampling trajectories that contribute to the estimator is minimal.

Typically, explicit solutions to MV-SDEs are unavailable. These SDEs are first approximated using a stochastic interacting particle system~\citep{mean_field_limit}, comprising $P$ coupled $d$-dimensional SDEs. This system is generally not solvable in closed form; hence, its trajectories are numerically approximated via a time-stepping scheme with $N$ time steps. Thus, the approximation of $\E{G(X(T))}$ for MV-SDEs involves two parameters: the numbers of particles $P$ and time steps $N$. Standard MC estimators based on the above $P$-particle approximation and the Euler--Maruyama time-stepping method~\citep{mlmc_mvsde,euler_mvsde} achieve $\order{\tol^{-4}}$ complexity for smooth, nonrare observables under bounded, Lipschitz coefficients. However, for rare events, the cost prohibitively increases owing to the exploding coefficient of variation problem~\citep{is_general_ref}. The {coefficient of variation} of the estimator explodes as the event becomes rarer, making standard MC intractable. To overcome this problem, importance sampling (IS) via a measure change in the path space is applied as a standard variance reduction technique~\citep{is_general_ref}. This technique involves transforming the underlying probability measure to sample rare trajectories and appropriately reweight them to ensure an unbiased estimator.

The initial studies of IS for MV-SDEs have appeared in the works by~\citep{is_mvsde} and~\citep{my_dlmcis}. A naive application of IS to the MV-SDE fails because, under a new measure, the process's law is no longer consistent with the original dynamics of the MV-SDE. Hence,~\citep{is_mvsde} proposed a decoupling approach, where the  law in the drift/diffusion of the MV-SDE is replaced by an empirical realization of the law precomputed from a particle system, after which a measure change is applied. Building on this idea, \citep{my_dlmcis} introduced a double-loop MC (DLMC) estimator and employed stochastic optimal control to design a time- and path-dependent IS control, reaching a complexity of $\order{\tol_\mathrm{r}^{-4}}$, matching the standard MC method for nonrare observables~\citep{mlmc_mvsde}. Recently, the multilevel MC method was combined with this IS scheme~\citep{my_mldlmcis}, reducing the complexity further compared with the DLMC estimator.

The multi-index MC (MIMC) method, introduced by~\citep{mimc2016}, extends the multilevel MC method by exploiting mixed refinement in multiple discretization directions (in this work, particles $P$ and time steps $N$). The MIMC method relies on mixed regularity with respect to these directions, further reducing the computational cost. In the MV-SDE setting, as $P \rightarrow \infty$, the errors from finite particle approximation and time discretization become coupled, an effect of the propagation of chaos~\citep{mean_field_limit} that cancels out leading-order error terms, which is leveraged in the multi-index complexity analysis. Although the MIMC method has been applied to smooth, nonrare observables of particle systems~\citep{mlmc_mvsde}, the goal of this  work is to combine the MIMC method with the IS scheme from~\citep{my_dlmcis}. This work presents an estimator developed for rare-event expectations of MV-SDEs that surpasses the multilevel DLMC estimator~\citep{my_mldlmcis}. Although the multilevel MC method with IS has been studied in other contexts (e.g., standard SDEs~\citep{mlmcis_kebaier}, ergodic SDEs~\citep{mlmcis_giles}, stochastic reaction networks~\citep{mlmcis_srn}, and diffusions in finance~\citep{mlmcis_alaya}), this work is the first integration of IS with the MIMC method to the best of our knowledge. The primary contributions of this work are as follows.
\begin{enumerate}
	\item This work extends the DLMC estimator of~\citep{my_dlmcis} to the multi-index setting, proposing a multi-index DLMC estimator for MV-SDEs. Until now, variance reduction for rare events in MV-SDEs had only been explored in single~\citep{is_mvsde,my_dlmcis} and multilevel MC~\citep{my_mldlmcis} settings. This work integrates IS with MIMC for MV-SDEs. This work applies the IS control from~\citep{my_dlmcis} across all multi-indices and introduces a strongly coupled antithetic sampler to control the rare-event variance of  the mixed-difference estimators. 
	\item This work derives a complexity result (Theorem~\ref{thm:midlmc_complexity}) for MV-SDEs based on assumptions on the bias and variance of the estimator (see Assumptions~\ref{ass:midlmc_bias}-\ref{ass:midlmc_var}) that are similar to standard mixed-difference assumptions made in MIMC. This work reveals that the proposed approach improves on the multilevel DLMC estimator, which achieves $\order{\tol_\mathrm{r}^{-3}}$ complexity, to an asymptotic complexity of $\order{\tol_{\mathrm{r}}^{-2} (\log \tol_{\mathrm{r}}^{-1})^2}$. This work presents strong numerical evidence in support of these assumptions on a simple MV-SDE, leaving further theoretical analysis to future work.
	\item This work numerically investigates the influence of IS on variance reduction in the proposed multi-index DLMC estimator. For sufficiently smooth observables, the numerical simulations demonstrate a significant variance reduction on all mixed-difference estimators. In the Kuramoto test case (Section~\ref{sec:kuramoto}), this result improves the complexity from $\order{\tol_{\mathrm{r}}^{-4}}$ obtained in the work by~\citep{my_dlmcis} to $\order{\tol_{\mathrm{r}}^{-2} (\log \tol_{\mathrm{r}}^{-1})^2}$ in the multi-index setting (see Figure~\ref{fig:adaptive_midlmcis}). This substantial improvement makes estimating rare-event quantities up to a relative error tolerance $\tol_{\mathrm{r}}$ feasible in practice. Moreover, this work numerically studies the effect of the smoothness of the observable on the mixed-difference bias and variance convergence rates. The experiments confirm that the assumptions hold for a sufficiently smooth $G$ and demonstrate how the convergence degrades as $G$ becomes less regular.
\end{enumerate}
This paper is organised as follows. Section~\ref{sec:mvsde} introduces the MV-SDE, associated notation, and problem statement.  Next, Section~\ref{sec:dlmc} reviews the decoupling approach in the work by~\citep{is_mvsde}, recalls the IS control for the decoupled MV-SDE from~\citep{my_dlmcis}, and formulates a DLMC estimator with IS. Then,  Section~\ref{sec:midlmc} presents the novel integration of IS with the multi-index DLMC estimator and derives the complexity bounds, conditional on two key assumptions on mixed-difference bias and variance decay (Assumptions~\ref{ass:midlmc_bias} and \ref{ass:midlmc_var}). Section~\ref{sec:numerics} applies the method to the Kuramoto model from statistical physics to test the proposed approach and empirically verify Assumptions~\ref{ass:midlmc_bias} and~\ref{ass:midlmc_var} for problems of varying regularity. The study confirms the theoretical complexity rates for the multi-index DLMC estimator, demonstrating the effectiveness of the proposed method in rare-event estimation for sufficiently smooth problems.

\textbf{Notation:} {The notation $\langle \cdot,\cdot \rangle$ represents the Euclidean dot product on $\mathbb{R}^d$, and $\norm{\cdot}$ denotes the corresponding norm. For a scalar function, $\nabla$ and $\nabla^2$ denote the gradient and Hessian, respectively. The notation $\cdot: \cdot$ signifies the Frobenius inner product between two matrix-valued functions. We denote by $C^k(A,B)$ the space of functions $f:A \rightarrow B$ that are $k$-times continuously differentiable on $A$. When $B = \mathbb{R}$, we simply write $C^k(A)$.} {For non-negative quantities $a$ and $b$, $a \lesssim b$ indicates the existence of a constant $c > 0$ such that $a \leq cb$. Unless stated otherwise, $c$ is independent of the discretisation and tolerance parameters, but may depend on fixed problem data.}
\section{McKean--Vlasov stochastic differential equation}
\label{sec:mvsde}

This work considers a broad class of McKean--Vlasov equations arising from the mean-field limit of stochastic interacting particle systems with pairwise interaction kernels~\citep{mean_field_limit}, defined on the probability space $\{\Omega,\mathcal{F},\{\mathcal{F}_t\}_{t \geq 0},\mathbb{P}\}$, where $\mathcal{F}_t$ represents the filtration of a standard Wiener process $\{W(t):t \in [0,T]\}$. For measurable functions $b:\mathbb{R}^d \cross \mathbb{R} \longrightarrow \mathbb{R}^d$, $\sigma:\mathbb{R}^d \cross \mathbb{R} \longrightarrow \mathbb{R}^{d \cross d}$, $\kappa_1: \mathbb{R}^d \cross \mathbb{R}^d \longrightarrow \mathbb{R}$, and $\kappa_2: \mathbb{R}^d \cross \mathbb{R}^d \longrightarrow \mathbb{R}$, this work studies the following It\^o SDE for the McKean--Vlasov stochastic process $X: [0,T] \times \Omega \rightarrow \mathbb{R}^d$:

\begin{empheq}[left=\empheqlbrace, right = ,]{equation} 
    \label{eqn:mvsde}
    \begin{alignedat}{2}
        \dd X(t) &= b\left(X(t),\int_{\mathbb{R}^d} \kappa_1 (X(t),x) \mu_t(\dd x)\right) \dd t  \\
        &\qquad + \sigma \left(X(t),\int_{\mathbb{R}^d} \kappa_2 (X(t),x) \mu_t(\dd x)\right) \dd W(t), \quad t>0 \\
        X(0) &= x_0 \sim \mu_0 \in \mathcal{P}(\mathbb{R}^d) , 
    \end{alignedat}
\end{empheq}
where $W:[0,T] \cross \Omega \longrightarrow \mathbb{R}^d$ represents a standard $d$-dimensional Wiener process with mutually independent components; $\mu_t \in \mathcal{P}(\mathbb{R}^d)$ denotes the law of $X(t)$, where $\mathcal{P}(\mathbb{R}^d)$ indicates the space of probability measures in $\mathbb{R}^d$; and $x_0 \in \mathbb{R}^d$ indicates a random initial state with the distribution $\mu_0 \in \mathcal{P}(\mathbb{R}^d)$. The functions $b$ and $\sigma$ represent the drift and diffusion functions/coefficients, respectively.

{
\textbf{Standing assumptions}: Throughout, unless stated otherwise, $b$ and $\sigma$ are globally Lipschitz in both variables with at most linear growth, $\kappa_1$ and $\kappa_2$ are bounded and globally Lipschitz, $\mu_0$ has sufficiently high finite moments for the weak and strong error estimates invoked in this work, and the Wiener process and initial state in~\eqref{eqn:mvsde} are mutually independent. These assumptions ensure the existence and uniqueness of solutions to~\eqref{eqn:mvsde} ~\citep{mean_field_limit}. The well-posedness of~\eqref{eqn:mvsde} under substantially weaker assumptions is known~\citep{mvsde_existence,mvsde_weak_soln,Chaudru-de-Raynal:2020aa,Chaudru-de-Raynal:2022aa}, but such problems are not the focus of this work. Additional smoothness assumptions required for the convergence rates stated in Assumptions~\ref{ass:midlmc_bias} and~\ref{ass:midlmc_var} are discussed in Section~\ref{sec:midlmc}.
}

The Fokker--Planck equation associated with~\eqref{eqn:mvsde} is a nonlinear integro - differential partial differential equation (PDE) containing nonlocal interaction terms~\citep{mvsde_pde}. Numerically solving this PDE to obtain solutions to~\eqref{eqn:mvsde} with a prescribed relative accuracy is computationally infeasible. Instead, this work employs strong approximations that converge in $L^p$ to the solution of~\eqref{eqn:mvsde}.

A strong approximation is given by a system of $P$ exchangeable It\^o SDEs, also known as a stochastic interacting particle system, with pairwise interaction kernels~\citep{mean_field_limit}. For $p=1, \ldots, P$, the state $X^P_p:[0,T]\cross\Omega \rightarrow \mathbb{R}^d$ of particle $p$ (out of $P$) follows the following SDE:

\begin{empheq}[left=\empheqlbrace, right = ,]{alignat=2}
    \label{eqn:strong_approx_mvsde}
    \dd X^P_p(t) &= b\left(X^P_p(t), \frac{1}{P} \sum_{j=1}^P  \kappa_1(X^P_p(t),X^P_j(t)) \right) \dd t \nonumber \\
    &\qquad + \sigma\left(X^P_p(t), \frac{1}{P} \sum_{j=1}^P \kappa_2(X^P_p(t),X^P_j(t)) \right) \dd W_p(t), \quad t>0 \\
    X^P_p(0) &= (x_0)_p \sim \mu_0 \in \mathcal{P}(\mathbb{R}^d) , \nonumber
\end{empheq}
where $\{(x_0)_p\}_{p=1}^P$ denotes independent and identically distributed (i.i.d.) random variables sampled from the initial distribution $\mu_0$, and $\{W_p\}_{p=1}^P$ represents mutually independent $d$-dimensional Wiener processes that are independent of $\{(x_0)_p\}_{p=1}^P$. Equation~\eqref{eqn:strong_approx_mvsde} approximates the law $\mu_t$ in \eqref{eqn:mvsde} using an empirical law based on $P$ particles $\{X^P_p(t)\}_{p=1}^P$:

\begin{equation}
    \label{eqn:emp_dist_law}
    \mu_t(\dd x) \approx \mu_t^P(\dd x) = \frac{1}{P} \sum_{j=1}^P \delta_{X^P_j(t)} (\dd x) \cdot
\end{equation}
In the limit $P \rightarrow \infty$, $\mu_t^P$ converges to $\mu_t$ (propagation of chaos) under suitable conditions~\citep{mean_field_limit}. Strong convergence (in the $L^p$-sense) of the particle system~\eqref{eqn:strong_approx_mvsde} to the MV-SDE~\eqref{eqn:mvsde} has been well established~\citep{mvsde_strong_conv_3,mvsde_strong_conv_2,mvsde_strong_conv_1}. The high dimensionality of the Fokker--Planck equation, satisfied by the joint probability density of the particle system, motivates the use of MC methods, which do not suffer from the curse of dimensionality.

\subsection{Motivating example: Kuramoto oscillator model}
\label{sec:kuramoto}

The method in this work is assessed on the simple, one-dimensional (1D) Kuramoto model \citep{chemical_app,neuroscience_app,combustion_app}, describing synchronization in statistical physics to model the behavior of large sets of fully connected, synchronized oscillators. The phase of oscillator $p$ (out of $P$) is represented by the process $X_p^P:[0, T]\cross\Omega \rightarrow \mathbb{R}$ with the following Itô SDE dynamics: 

\begin{empheq}[left=\empheqlbrace, right = ,]{alignat=2}
    \label{eqn:kuramoto_model}
    \dd X^P_p(t) &= \left(\xi_p + \frac{1}{P} \sum_{q=1}^P \sin\left(X^P_p(t) - X^P_q(t)\right)\right) \dd t + \sigma \dd W_p (t) , \quad t>0\\
    X^P_p(0) &= (x_0)_p \sim \mu_0 \in \mathcal{P}(\mathbb{R}) , \nonumber
\end{empheq}
where $\{\xi_p\}_{p=1}^P$ denotes i.i.d. random variables sampled from a prescribed distribution and models the natural frequency of each oscillator. The diffusion $\sigma \in \mathbb{R}$ is constant, and $\{(x_0)_p\}_{p=1}^P$ represents i.i.d. random variables sampled from $\mu_0 \in \mathcal{P}(\mathbb{R})$. In addition, $\{W_p\}_{p=1}^P$ represents mutually independent 1D Wiener processes, and $\{\xi_p\}_{p=1}^P, \{(x_0)_p\}_{p=1}^P$, $\{W_p\}_{p=1}^P$ are mutually independent. The coupled particle system~\eqref{eqn:kuramoto_model} strongly converges in the mean-field limit as $P \rightarrow \infty$, where each particle satisfies the following MV-SDE:

\begin{empheq}[left=\empheqlbrace, right = ,]{equation}
    \label{eqn:kuramoto_mvsde}
    \begin{alignedat}{2}
        \dd X(t) &= \left(\xi + \int_{\mathbb{R}} \sin (X(t)-x) \mu_t (\dd x) \right) \dd t + \sigma \dd W(t), \quad t>0 \\
        X(0) &= x_0 \sim \mu_0 \in \mathcal{P}(\mathbb{R}) , 
    \end{alignedat}
\end{empheq}
where $X(t)$ denotes the McKean--Vlasov process at time $t$, $\xi$ represents a random variable sampled from a prescribed distribution, and $\mu_t$ indicates the law of $X(t)$. In the notation of~\eqref{eqn:mvsde}, this example corresponds to $b(x,y) = \xi + y$, $\sigma(x,y) = \sigma$, $\kappa_1(x,y) = \sin(x-y)$ and $\kappa_2(x,y) = 0$, satisfying the standard conditions ensuring the well-posedness of~\eqref{eqn:kuramoto_mvsde}. This model is used as a test-bed because it is a simple, well-understood 1D MV-SDE that can be tuned by adjusting the parameters, such as the coupling kernel, and the distribution of $\xi$, to control the rarity of the event of interest.

\subsection{Problem setting}
\label{sec:obj}

Let $T>0$ be a finite terminal time and $X:[0,T] \cross \Omega \rightarrow \mathbb{R}^d$ denote the McKean--Vlasov process~\eqref{eqn:mvsde}. Let $G: \mathbb{R}^d \longrightarrow \mathbb{R}$ be a given nonnegative scalar observable that satisfies $G(x) \geq 0$ for all $x \in \mathbb{R}^d$. The objective is to build a computationally efficient MC estimator $\mathcal{A}$ of $\E{G(X(T))}$ that satisfies a relative error tolerance $\tol_{\mathrm{r}} > 0$ for a given confidence level determined by $0 < \nu \ll 1$ in the following sense:
\begin{equation}
    \label{eqn:objective}
    \prob{\frac{\abs{\mathcal{A}-\E{G(X(T))}}}{\abs{\E{G(X(T))}}} \geq \tol_{\mathrm{r}}} \leq \nu \cdot
\end{equation}
This work aims to compute rare events associated with the final state $X(T)$, implying that the distribution of $G(X(T))$ places most of its probability mass deep in the tail of the law of $X(T)$. The computational feasibility of the vanilla MC method rapidly diminishes due to the exploding coefficient of variation~\citep{is_general_ref}. Hence, the MC method must be combined with the IS variance reduction technique to produce feasible estimates of rare events. However, applying this technique to MV-SDEs raises a considerable difficulty. The drift and diffusion functions are linked to $\{\mu_t: t \in [0,T]\}$, the law of the McKean--Vlasov process. Under any IS measure change (e.g., via Girsanov's theorem), $\mu_t$ is no longer the law of the modified process. This challenge is overcome using the decoupling approach introduced by~\citep{is_mvsde} by transforming the MV-SDE~\eqref{eqn:mvsde} into a standard SDE by conditioning it on an empirical approximation of $\mu_t$. In summary, computing rare-event quantities for MV-SDEs requires variance reduction (to control the coefficient of variation explosion) and careful handling of mean-field coupling under measure change.  Section~\ref{sec:dlmc} introduces the decoupling approach for MV-SDEs that addresses these problems.

\section{Double-loop importance sampling for McKean--Vlasov stochastic differential equations}
\label{sec:dlmc}

The decoupling approach replaces the deterministic law $\{\mu_t: t \in [0,T]\}$ in~\eqref{eqn:mvsde} with an empirical approximation~\eqref{eqn:emp_dist_law} using the interacting particle system~\eqref{eqn:strong_approx_mvsde}. Given this precomputed empirical law, a decoupled MV-SDE is introduced, and an IS measure change is applied to it. Hence, the empirical  estimation of the law $\{\mu_t: t \in [0,T]\}$ and the IS measure change are decoupled. The decoupling approach involves the following steps~\citep{my_dlmcis}:
\begin{enumerate}
	\item The law $\{\mu_t:t \in [0,T]\}$ in~\eqref{eqn:mvsde} is approximated using the empirical measure $\{\mu_t^P: t \in [0,T]\}$ in~\eqref{eqn:emp_dist_law} using one realization of the $P$-particle system \eqref{eqn:strong_approx_mvsde} with particles $\{X^P_p\}_{p=1}^P$.
	\item Given $\{\{X^P_p(t)\}_{p=1}^P \sim \mu_t^P: t \in [0,T]\}$, this work defines the decoupled McKean--Vlasov process $\Bar{X}^P:[0,T] \times \Omega \rightarrow \mathbb{R}^d$ as the solution to the following It{\^o} SDE:
	\begin{empheq}[left=\empheqlbrace, right = ,]{alignat=2}
        \label{eqn:decoupled_mvsde}
        \begin{split}
            &\dd \Bar{X}^P(t) = b\left(\Bar{X}^P(t), \frac{1}{P} \sum_{j=1}^P  \kappa_1(\Bar{X}^P(t),X^P_j(t)) \right) \dd t \\
            &\qquad + \sigma \left(\Bar{X}^P(t), \frac{1}{P} \sum_{j=1}^P  \kappa_2(\Bar{X}^P(t),X^P_j(t)) \right) \dd \Bar{W}(t), \quad t \in [0,T] \\
            &\Bar{X}^P(0) = \Bar{x}_0 \sim \mu_0, \quad \Bar{x}_0 \in \mathbb{R}^d , 
        \end{split}
    \end{empheq}
    where the superscript $P$ indicates that the drift and diffusion functions in~\eqref{eqn:decoupled_mvsde} are computed using $\{ \mu_t^P: t \in [0,T]\}$ derived from the stochastic $P$-particle system $\{X^P_p\}_{p=1}^P$. In addition, $\Bar{W}$ denotes a standard $d$-dimensional Wiener process independent of the Wiener processes $\{W_p\}_{p=1}^P$ employed in~\eqref{eqn:strong_approx_mvsde}, and $\Bar{x}_0 \in \mathbb{R}^d$ represents a random initial state sampled from $\mu_0$, independent of $\{(x_0)_p\}_{p=1}^P$ in~\eqref{eqn:strong_approx_mvsde}. Thus, \eqref{eqn:decoupled_mvsde} is a standard SDE for a given realization of the empirical law $\{\mu_t^P: t \in [0, T]\}$.
    \item Let the $P$-particle system~\eqref{eqn:strong_approx_mvsde} $\{X^P_p\}_{p=1}^P : [0,T] \times \Omega \rightarrow (\mathbb{R}^d)^P$ be defined on the probability space $(\Omega,\mathcal{F},\mathbb{P})$ with $\mathbb{E}_\mathbb{P}$ denoting the expectation with respect to $\mathbb{P}$. For given $\{X^P_p\}_{p=1}^P$, this work defines $\bar{X}^P$ in~\eqref{eqn:decoupled_mvsde} on the copy space $(\tilde{\Omega},\tilde{\mathcal{F}},\tilde{\mathbb{P}})$. Hence, $\bar{X}^P$ is defined on the product space $(\Omega,\mathcal{F},\mathbb{P}) \times (\tilde{\Omega},\tilde{\mathcal{F}},\tilde{\mathbb{P}})$ with $\mathbb{E}_{\mathbb{P} \otimes \tilde{\mathbb{P}}}$ denoting the expectation with respect to the product measure space $\mathbb{P} \otimes \tilde{\mathbb{P}}$.
	\item The quantity of interest is approximated as the following nested expectation:
    \begin{align}
        \label{eqn:total_exp_mvsde}
        \E{G(X(T))} & \approx \mathbb{E}_{\mathbb{P} \otimes \tilde{\mathbb{P}}} \left[G(\Bar{X}^P(T))\right] \nonumber \\
        &= \mathbb{E}_\mathbb{P} \left[\mathbb{E}_{\tilde{\mathbb{P}}} \left[ G(\Bar{X}^P(T)) \mid \{\mu_t^P: t \in [0,T]\} \right]\right] \cdot
    \end{align}
    Henceforth, $\E{G(\Bar{X}^P(T))} \equiv \mathbb{E}_{\mathbb{P} \otimes \tilde{\mathbb{P}}} \left[G(\Bar{X}^P(T))\right]$ for convenience.
\end{enumerate}
As~\eqref{eqn:decoupled_mvsde} is a standard SDE for a given empirical law $\{\mu_t^P: t \in [0,T]\}$, this work follows the procedure from~\citep{my_dlmcis} to derive the optimal IS measure change.

\subsection{Importance sampling using stochastic optimal control for decoupled McKean--Vlasov stochastic differential equation}
\label{sec:soc_mvsde}

The PDE that yields the optimal control for the decoupled MV-SDE for a given empirical law $\{\mu_t^P: t \in [0,T]\}$ was formulated by~\citep{my_dlmcis}. This work defines the controlled process $\Bar{X}^P_\zeta: [0,T] \cross \Omega \rightarrow \mathbb{R}^d$ following the given controlled dynamics with the Markov control function $\zeta:[0,T] \cross \mathbb{R}^d \rightarrow \mathbb{R}^d$:
    \begin{empheq}[left=\empheqlbrace, right = \cdot]{alignat=2}
        \label{eqn:dmvsde_sde_is}
        \begin{split}
            \dd \Bar{X}^P_\zeta(t) &= \Bigg( \Bigg. b\Bigg(\Bar{X}^P_\zeta(t), \frac{1}{P} \sum_{j=1}^P  \kappa_1(\Bar{X}^P_\zeta(t),X^P_j(t)) \Bigg) \\
            &\qquad + \sigma \Bigg(\Bar{X}^P_\zeta(t), \frac{1}{P} \sum_{j=1}^P  \kappa_2(\Bar{X}^P_\zeta(t),X^P_j(t)) \Bigg) \zeta(t,\Bar{X}^P_\zeta(t)) \Bigg. \Bigg) \dd t \\
            &\quad + \sigma \left(\Bar{X}^P_\zeta(t), \frac{1}{P} \sum_{j=1}^P  \kappa_2(\Bar{X}^P_\zeta(t),X^P_j(t)) \right) \dd W(t),  0<t<T \\ 
            \Bar{X}^P_\zeta(0) &= \Bar{X}^P(0) = \Bar{x}_0 \sim \mu_0 . 
        \end{split}
    \end{empheq}
This work applies~\eqref{eqn:strong_approx_mvsde} to obtain $\{\{X^P_p(t)\}_{p=1}^P \sim \mu^P_t:t \in [0,T]\}$ in~\eqref{eqn:dmvsde_sde_is}. This work assumes that $\zeta$ satisfies Novikov's condition:
    \begin{equation}
        \label{eqn:novikov_cond}
    	\E{\exp{\frac{1}{2} \int_0^T \norm{\zeta(s,\bar{X}^P_\zeta(s))}^2 \dd s} \quad \Bigg| \quad \{\mu^P_t:t \in [0,T]\}} < \infty \cdot
    \end{equation}  
The value function $u:[0,T] \times \mathbb{R}^d \rightarrow \mathbb{R}$ that minimizes the second moment of the MC estimator of  $\mathbb{E}\left[G(\Bar{X}^P(T)) \mid \{\mu_t^P: t \in [0, T]\}\right]$ is defined as follows:
    \begin{align}
        \label{eqn:dmvsde_value_fxn}
        u(t,x) &= \min_{\zeta \in \mathcal{Z}} \mathbb{E}\Bigg[\Bigg.G^2(\Bar{X}_\zeta^P(T)) \exp \Bigg\{ \Bigg. -\int_t^T \norm{\zeta(s,\Bar{X}_\zeta^P(s))}^2 \dd s \nonumber\\
        &\qquad \quad - 2 \int_t^T \langle \zeta(s,\Bar{X}_\zeta^P(s)),\dd W(s) \rangle \Bigg. \Bigg\} \Bigg| \Bar{X}_\zeta^P(t) = x, \{\mu_t^P: t \in [0,T]\}\Bigg.\Bigg] \cdot
    \end{align} 
{$\mathcal{Z}$ represents the class of deterministic Markov feedback functions $\zeta:[0,T] \cross \mathbb{R}^d \rightarrow \mathbb{R}^d$ that are sufficiently regular to ensure well-posedness of~\eqref{eqn:dmvsde_sde_is}, satisfy Novikov's condition~\eqref{eqn:novikov_cond} conditional on the empirical law, and make the second moment in~\eqref{eqn:dmvsde_value_fxn} finite.} This work defines {$v:[0,T] \cross \mathbb{R}^d \rightarrow \mathbb{R}$} such that $u(t, x) = v^2(t, x)$. Then, $v$ satisfies the following linear backward PDE: 
    \begin{empheq}[left=\empheqlbrace, right = ,]{alignat=2}
	    \label{eqn:dmvsde_hjb_form3}
	    \begin{split}
	        &\frac{\partial v}{\partial t} + \langle b\left(x, \frac{1}{P} \sum_{j=1}^P  \kappa_1(x,X^P_j(t)) \right), \nabla v \rangle \\
	        &+ \frac{1}{2} \nabla^2 v : \left(\sigma \sigma^T\right) \left(x, \frac{1}{P} \sum_{j=1}^P  \kappa_2(x,X^P_j(t)) \right) = 0 , (t,x) \in [0,T) \cross \mathbb{R}^d \\
	        &v(T,x) = \abs{G(x)}, \quad x \in \mathbb{R}^d , 
	    \end{split}
	\end{empheq}
    with the following optimal control:
    \begin{equation}
	    \label{eqn:dmvsde_hjb_optimal_control}
	    \zeta^*(t,x) = \sigma^T \left(x, \frac{1}{P} \sum_{j=1}^P  \kappa_2(x,X^P_j(t)) \right) \nabla \log v (t,x),
	\end{equation}
which minimizes the second moment in~\eqref{eqn:dmvsde_value_fxn}, for a given $\{\{X^P_p(t)\}_{p=1}^P \sim \mu_t^P: t \in [0,T]\}$. The corresponding proof is available in Appendix B in the work by~\citep{my_dlmcis}. {When $G$ vanishes for the indicator function at $t=T$, $\nabla \log v$ is interpreted only on $[0,T)$ where the parabolic smoothing of~\eqref{eqn:dmvsde_hjb_form3} yields $v > 0$ under non-degeneracy assumptions. In numerical experiments, we evaluate the control $\zeta^*$ only at time nodes less than $T$.} The IS control $\zeta^*$, which depends on the time and state, tilts the paths of~\eqref{eqn:dmvsde_sde_is} towards the rare-event domain and reduces the variance of the MC estimator. The likelihood resulting from the IS measure change is written as the following Radon--Nikodym derivative:
\begin{equation}
	\label{eqn:radon_nikodym}
	\mathbb{L}_{[0,T]}^P = \exp{-\int_0^T \langle \zeta(s, \bar{X}_\zeta^P(s)), \dd W(s) \rangle - \frac{1}{2} \int_0^T \norm{\zeta(s, \bar{X}_\zeta^P(s))}^2 \dd s} \cdot
\end{equation}
With this measure change, the nested expectation~\eqref{eqn:total_exp_mvsde} is rewritten as follows:
\begin{align}
	\E{G(X(T))} & \approx \mathbb{E}_\mathbb{P} \left[\mathbb{E}_{\tilde{\mathbb{P}}} \left[ G(\Bar{X}^P(T)) \mid \{\mu_t^P: t \in [0,T]\} \right]\right] \nonumber\\
	\label{eqn:total_exp_mvsde_is}
	&= \mathbb{E}_\mathbb{P} \left[\mathbb{E}_{\tilde{\mathbb{Q}}} \left[ G(\Bar{X}^P_\zeta(T)) \mathbb{L}_{[0,T]}^P \mid \{\mu_t^P: t \in [0,T]\} \right]\right],
\end{align}
where $\tilde{\mathbb{Q}}$ denotes the shifted probability measure of the copy space, such that~\eqref{eqn:dmvsde_sde_is} is defined on $(\Omega,\mathcal{F},\mathbb{P}) \cross (\tilde{\Omega},\tilde{\mathcal{F}},\tilde{\mathbb{Q}})$. The nested expectation~\eqref{eqn:total_exp_mvsde_is} is approximated using a DLMC estimator.

\subsection{Double loop Monte Carlo estimator with importance sampling}
\label{sec:dlmcis}

The DLMC estimator for the nested expectation~\eqref{eqn:total_exp_mvsde_is} is constructed as follows.

\begin{enumerate}
	\item The empirical law $\{\mu_t^P: t \in [0,T]\}$ approximates $\mu_t$ in~\eqref{eqn:mvsde} using~\eqref{eqn:strong_approx_mvsde}. In practice, a time-discretized empirical law with $N$ time steps is obtained from the Euler--Maruyama discretization $\{X_p^{P \mid N}\}_{p=1}^P$ of the particle system. The Euler--Maruyama continuous-time extension of the time-discretized stochastic particle system $\{X_p^{P \mid N}\}_{p=1}^P$ is applied to extend the time-discrete empirical law to all $t \in [0,T]$. Then, the empirical law is defined as follows:
	\begin{equation}
		\label{eqn:emp_law_em_cont_time}
		\mu^{P \mid N}(t) = \frac{1}{P} \sum_{j=1}^P \delta_{X_j^{P \mid N}(t)}, \quad \forall t \in [0,T] \cdot
	\end{equation}
	The $P$ underlying sets of random variables employed to generate a realization of $\{X_p^{P \mid N}\}_{p=1}^P$ using~\eqref{eqn:strong_approx_mvsde} are denoted by $\omega_{1:P}$.
	\item Given $\mu^{P \mid N}$ from~\eqref{eqn:emp_law_em_cont_time},~\eqref{eqn:dmvsde_hjb_form3} is solved and~\eqref{eqn:dmvsde_hjb_optimal_control} is employed to obtain the optimal Markov control function $\zeta$ for IS.
	\item Given $\mu^{P \mid N}$ from~\eqref{eqn:emp_law_em_cont_time} and the control function $\zeta$, sample paths of the controlled decoupled MV-SDE~\eqref{eqn:dmvsde_sde_is} are generated. The Euler--Maruyama discretization of SDE~\eqref{eqn:dmvsde_sde_is} with $N$ time steps is applied to obtain the approximate sample paths denoted by $\{\bar{X}_\zeta^{P \mid N}(t_n)\}_{n=1}^N$. The set of random variables employed to generate a sample path is denoted by $\tilde{\omega}$.
	\item The DLMC estimator $\mathcal{A}_\mathrm{MC}$ approximates the nested expectation~\eqref{eqn:total_exp_mvsde_is} using the following nested sample average:
	\begin{equation}
        \label{eqn:dmvsde_dlmc_est}
        \mathcal{A}_{\textrm{MC}} = \frac{1}{M_1} \sum_{i=1}^{M_1} \frac{1}{M_2} \sum_{j=1}^{M_2} G \left(\Bar{X}_\zeta^{P|N}(T) \right) \mathbb{L}^{P|N} \left( \omega_{1:P}^{(i)} \cross {\tilde{\omega}^{(i,j)}} \right) ,
    \end{equation}
    where $\mathbb{L}^{P|N}$ denotes the time-discretized approximation of the likelihood~\eqref{eqn:radon_nikodym} defined as follows:
    \begin{align}
        \label{eqn:dlmc_llhood_factor}
        \mathbb{L}^{P|N} &= \prod_{n=0}^{N-1} \exp \Bigg\{ \Bigg. -\frac{1}{2} \Delta t \norm{\zeta(t_n,\Bar{X}_\zeta^{P|N}(t_n))}^2 \\
        &\qquad - \langle \Delta W(t_n), \zeta(t_n,\Bar{X}^{P|N}_\zeta(t_n)) \rangle \Bigg. \Bigg\} \nonumber,
    \end{align}
     where $M_1$ denotes the number of independent realizations of $\mu^{P|N}$ and $\omega_{1:P}^{(i)}$ denotes the $i^{\mathrm{th}}$ realization of $\omega_{1:P}$. For each realization of $\mu^{P|N}$, $M_2$ represents the number of sample paths of the decoupled MV-SDE, and ${\tilde{\omega}^{(i,j)}}$ denotes the $j^{\mathrm{th}}$ realization of $\tilde{\omega}$ given $\omega_{1:P}^{(i)}$. Further, $\Delta W(t_n)$ represents the Wiener increment driving the dynamics of the time-discretized decoupled MV-SDE~\eqref{eqn:dmvsde_sde_is}.
\end{enumerate}

Algorithm~\ref{alg:dlmc_outline} in~\ref{app:dlmc_outline} summarises the above approach as pseudocode. {In the ideal double loop formulation, the PDE~\eqref{eqn:dmvsde_hjb_form3} should be numerically solved for each realisation of $\mu^{P \mid N}$. In Algorithm~\ref{alg:dlmc_outline}, we instead compute a fixed off-line control $\zeta$ from one sufficiently fine empirical law $\mu^{\bar{P} \mid \bar{N}}$ and reuse the same admissible control for all outer loop samples. This does not affect the bias of the IS estimator for all discretised decoupled MV-SDE~\eqref{eqn:dmvsde_sde_is} simulations, but it makes the control suboptimal relative to the realisation-dependent optimiser in~\eqref{eqn:dmvsde_hjb_optimal_control}. The quality of this approximation is assessed numerically in Section~\ref{sec:numerics}.} Under standard assumptions, this DLMC estimator with IS~\eqref{eqn:dmvsde_dlmc_est} achieves $\mathcal{O}(\tol_{\mathrm{r}}^{-4})$ complexity for a relative error tolerance $\tol_{\mathrm{r}}$~\citep{my_dlmcis}, which is a significant improvement over the standard MC method that is infeasible. Although the DLMC estimator mitigates the explosive growth of the coefficient of variation and achieves the same complexity order as the standard MC method for nonrare observables, further improvements are possible. In particular, the DLMC method still samples all trajectories at a fixed high resolution ($P,N$). In the work by~\citep{my_mldlmcis}, the DLMC estimator was extended to the multilevel MC setting. This extension allows for  estimating the expectations of Lipschitz rare-event observables up to a relative tolerance of $\tol_{\mathrm{r}}$ while reducing the complexity to $\mathcal{O}(\tol_{\mathrm{r}}^{-3})$ for the Kuramoto example \eqref{eqn:kuramoto_model}. The presence of two discretization parameters, $P$ and $N$, in MV-SDEs motivates the extension of this estimator to an MIMC setting to further reduce the complexity.

\section{Multi-index double loop Monte Carlo method}
\label{sec:midlmc}

This section aims to combine the multi-index technique with the stochastic optimal control-based IS scheme for rare events introduced in Section~\ref{sec:soc_mvsde}. This approach yields the novel multi-index DLMC estimator. Following the work by~\citep{mimc2016}, this work introduces the following multi-index DLMC discretization. Two discretization parameters $(P, N)$ are required to generate sample paths of the decoupled MV-SDE~\eqref{eqn:dmvsde_sde_is}. The multi-index $\boldsymbol\alpha = \left( \alpha_1, \alpha_2 \right) \in \mathbb{N}^2$ defines the discretization parameters as follows:

\begin{equation}
\label{eqn:mimc_discr}
	P_{\alpha_1} = P_0 \times 2^{\alpha_1}, \quad N_{\alpha_2} = N_{0} \times 2^{\alpha_2} ,
\end{equation}
where $P_0$ and $N_0$ represent the minimum number of particles and time steps, respectively, to generate the approximate sample paths of the decoupled MV-SDE~\eqref{eqn:dmvsde_sde_is}. This work proposes the following method to couple IS with the MIMC estimator. The optimal control PDE~\eqref{eqn:dmvsde_hjb_form3} is solved using one realization of the stochastic particle system~\eqref{eqn:strong_approx_mvsde} with a large number of particles $\bar{P}$ and time steps $\bar{N}$, obtaining the approximate IS control function $\zeta$ using~\eqref{eqn:dmvsde_hjb_optimal_control}. The same policy $\zeta$ is applied to generate samples of the decoupled MV-SDE~\eqref{eqn:dmvsde_sde_is} at all discretization multi-indices $\boldsymbol\alpha$. Let $G_{\boldsymbol\alpha} = G(\bar{X}_\zeta^{P_{\alpha_1}|N_{\alpha_2}}(T)) \mathbb{L}^{P_{\alpha_1}|N_{\alpha_2}}$, where $\bar{X}_\zeta^{P_{\alpha_1}|N_{\alpha_2}}$ denotes the Euler--Maruyama discretized controlled decoupled MV-SDE process~\eqref{eqn:dmvsde_sde_is} conditioned on the empirical law $\mu^{P_{\alpha_1}|N_{\alpha_2}}$~\eqref{eqn:emp_law_em_cont_time}. The discretized likelihood factor $\mathbb{L}^{P_{\alpha_1}|N_{\alpha_2}}$ is computed using~\eqref{eqn:dlmc_llhood_factor}. Following the work by~\citep{mimc2016}, this work defines the first-order mixed-difference operator in two dimensions for this setting: 

\begin{equation}
	\label{eqn:mixed_diff}
	{\Delta^\text{mix} G_{\boldsymbol\alpha}} = \left( G_{(\alpha_1,\alpha_2)} - G_{(\alpha_1-1,\alpha_2)} \right) - \left( G_{(\alpha_1,\alpha_2-1)} - G_{(\alpha_1-1,\alpha_2-1)} \right) \cdot
\end{equation}
{Assuming that $G_{\boldsymbol\alpha}$ converges to $G(X(T))$  (in $L^1$-sense) as $\min(\alpha_1,\alpha_2) \rightarrow \infty$, the multi-index telescoping identity gives
}
\begin{equation}
	\label{eqn:mimc_concept}
    {
	\E{G(X(T))} = \lim_{L_1, L_2 \rightarrow \infty} \sum_{\alpha_1=0}^{L_1} \sum_{\alpha_2=0}^{L_2} \E{{\Delta^\text{mix} G_{\boldsymbol\alpha}}} =  \sum_{\boldsymbol\alpha \in \mathbb{N}^2} \E{{\Delta^\text{mix} G_{\boldsymbol\alpha}}},
    }
\end{equation}
where $G_{(-1,0)} = 0$, $G_{(0,-1)} = 0$, and $G_{(-1,-1)} = 0$. Let ${\hat{\Delta} G_{\boldsymbol\alpha}}$ be a random variable such that $\E{{\hat{\Delta} G_{\boldsymbol\alpha}}}=\E{{\Delta^\text{mix} G_{\boldsymbol\alpha}}}$. In the trivial case, ${\hat{\Delta} G_{\boldsymbol\alpha}} = {\Delta^\text{mix} G_{\boldsymbol\alpha}}$. Alternatively, ${\hat{\Delta} G_{\boldsymbol\alpha}}$ can be chosen such that $\Var{{\hat{\Delta} G_{\boldsymbol\alpha}}} \ll \Var{{\Delta^\text{mix} G_{\boldsymbol\alpha}}}$. Each expectation in \eqref{eqn:mimc_concept} is  approximated using a DLMC estimator, resulting in the following multi-index DLMC estimator:

\begin{equation}
	\label{eqn:midlmc_estimator}
	\mathcal{A}_{\mathrm{MIMC}}(\mathcal{I}) = \sum_{\boldsymbol\alpha \in \mathcal{I}} \frac{1}{M_{1,\boldsymbol\alpha}} \sum_{m_1=1}^{M_{1,\boldsymbol\alpha}} \frac{1}{M_{2,\boldsymbol\alpha}} \sum_{m_2=1}^{M_{2,\boldsymbol\alpha}} {\hat{\Delta} G_{\boldsymbol\alpha}} \left( \omega_{1:P_{\alpha_1}}^{(\boldsymbol\alpha,m_1)} \times {\tilde{\omega}^{(\boldsymbol\alpha,m_1,m_2)}} \right) ,
\end{equation}
where $\mathcal{I} \subset \mathbb{N}^2$ represents an appropriately chosen index set, and $\{M_{1,\boldsymbol\alpha}, M_{2,\boldsymbol\alpha}\}$ denotes the integer numbers of samples in the inner and outer loops, respectively, of the DLMC estimator for each $\boldsymbol\alpha \in \mathcal{I}$. {Assumption~\ref{ass:midlmc_bias} is used to control the truncation bias for the admissible index set $\mathcal{I}$.} To control the variance of the mixed-differences, $G_{\boldsymbol\alpha}$ must be tightly coupled with the three lower-index terms in~\eqref{eqn:mixed_diff}. Based on the antithetic estimators in the studies by~\citep{mlmc_mvsde}, \citep{Szpruch:2021aa}, and~\citep{my_mldlmcis}, this coupling is achieved via an antithetic estimator ${\hat{\Delta} G_{\boldsymbol\alpha}}$ that reuses a large portion of the random variables (particles and Wiener paths) between multi-indices $\boldsymbol\alpha$ and their neighbors. The estimator ${\hat{\Delta} G_{\boldsymbol\alpha}}$ is defined as follows:
\begin{align}
	\label{eqn:antithetic_sampler}
	{\hat{\Delta} G_{\boldsymbol\alpha}} \left( \omega_{1:P_{\alpha_1}}^{(\boldsymbol\alpha,m_1)} \times {\tilde{\omega}^{(\boldsymbol\alpha,m_1,m_2)}} \right) &= \Bigg( \Bigg. \left( G_{(\alpha_1,\alpha_2)} - \mathcal{G}_{(\alpha_1-1,\alpha_2)} \right) \\
	&\hspace{-5em}- \left( G_{(\alpha_1,\alpha_2-1)} - \mathcal{G}{(\alpha_1-1,\alpha_2-1)} \right) \Bigg. \Bigg) \left( \omega_{1:P_{\alpha_1}}^{(\boldsymbol\alpha,m_1)} \times {\tilde{\omega}^{(\boldsymbol\alpha,m_1,m_2)}} \right), \nonumber
\end{align}
where $\mathcal{G}_{(\alpha_1-1,\alpha_2)}$ is highly correlated with $G_{(\alpha_1,\alpha_2)}$, which is defined as follows:

\begin{align}
	\label{eqn:antithetic_coupling}
	\mathcal{G}_{(\alpha_1-1,\alpha_2)} \left( \omega_{1:P_{\alpha_1}}^{(\boldsymbol\alpha,m_1)} \times {\tilde{\omega}^{(\boldsymbol\alpha,m_1,m_2)}} \right) &= \frac{1}{2} \left( G_{(\alpha_1-1,\alpha_2)} \Bigg( \Bigg. \omega_{1:P_{\alpha_1-1}}^{(\boldsymbol\alpha,m_1)} \times {\tilde{\omega}^{(\boldsymbol\alpha,m_1,m_2)}} \right) \\
	&+ G_{(\alpha_1-1,\alpha_2)} \left( \omega_{P_{\alpha_1-1}+1:P_{\alpha_1}}^{(\boldsymbol\alpha,m_1)} \times {\tilde{\omega}^{(\boldsymbol\alpha,m_1,m_2)}} \right) \Bigg. \Bigg) \cdot \nonumber
\end{align}
In \eqref{eqn:antithetic_coupling}, the $P_{\alpha_1}$ sets of random variables are split into two i.i.d. groups of size $P_{\alpha_1-1}$ and each group is employed to generate a realization of the empirical law independently. Given each realization of the empirical law, approximate sample paths of $\bar{X}_\zeta^{P_{\alpha_1-1} \mid N_{\alpha_2}}$ are generated using the same $\tilde{\omega}$ as for $G_{(\alpha_1,\alpha_2)}$ before averaging the quantity of interest over the two groups. Intuitively, six simulations in~\eqref{eqn:antithetic_sampler} are coupled at different discretizations in a carefully arranged way to yield a small variance.

\subsection{Error analysis}

The relative error of the multi-index DLMC estimator $\mathcal{A}_{\mathrm{MIMC}}$ is bounded as follows:
\begin{equation}
	\label{eqn:rel_error_split}
	\frac{\abs{\E{G(X(T))} - \mathcal{A}_{\mathrm{MIMC}}}}{\abs{\E{G(X(T))}}} \leq \underbrace{\frac{\abs{\E{G(X(T))} - \E{\mathcal{A}_{\mathrm{MIMC}}}}}{\abs{\E{G(X(T))}}}}_{=\epsilon_b,\text{ Relative bias}} + \underbrace{\frac{\abs{\E{\mathcal{A}_{\mathrm{MIMC}}} - \mathcal{A}_{\mathrm{MIMC}}}}{\abs{\E{G(X(T))}}}}_{=\epsilon_s,\text{ Relative statistical error}} \cdot
\end{equation}
Instead of satisfying a prescribed relative error tolerance $\tol_\mathrm{r}$ in the sense of~\eqref{eqn:objective}, the accuracy is split between the relative bias and statistical errors using parameter $\theta \in (0,1)$ and imposing the following stricter constraints for a confidence level determined by $0 < \nu \ll 1$:

\begin{align}
	\label{eqn:midlmc_bias_constraint}
	\epsilon_b \leq (1-\theta) \tol_{\mathrm{r}}, \\
	\label{eqn:midlmc_stat_constraint}
	\prob{\epsilon_s \geq \theta \tol_{\mathrm{r}}} \leq \nu \cdot
\end{align}
Throughout this work, $\theta$ is assumed to be given and fixed. This work makes the following assumptions regarding the first-order mixed-differences.
\begin{assump}{ (\textit{Multiplicative bias decay})}
\label{ass:midlmc_bias}
	\hspace{1mm} For all $\boldsymbol\alpha \in \mathcal{I}$, for constants $Q_B > 0$ and $b_1,b_2 > 0$, the absolute value of the expected value of ${\hat{\Delta} G_{\boldsymbol\alpha}}$ satisfies the following:
	\begin{equation}
	\label{eqn:midlmc_bias_ass}
		\abs{\E{{\hat{\Delta} G_{\boldsymbol\alpha}}}} \leq Q_B P_{\alpha_1}^{-b_1} N_{\alpha_2}^{-b_2} \cdot
	\end{equation}
\end{assump}
Assumption~\ref{ass:midlmc_bias} is motivated by the standard weak error expansion of the Euler--Maruyama time-discretization scheme~\citep{talay1990expansion} with respect to the number of time steps $N$, and the weak convergence of the interacting particle system approximation~\eqref{eqn:strong_approx_mvsde} with respect to the number of particles $P$~\citep{Mischler:2015aa,Bencheikh:2019aa,Szpruch:2021aa,Chassagneux:2022aa}. Based on~\citep{talay1990expansion}, $b_2=1$ is expected for the Euler--Maruyama scheme. A recent analysis by~\citep{theory_mvsde} established a $\order{P^{-1}}$ weak convergence rate of the particle system~\eqref{eqn:strong_approx_mvsde} to the MV-SDE~\eqref{eqn:mvsde}, indicating that $b_1 = 1$ under stronger conditions on $b,\sigma,\kappa_1,\kappa_2,G$ than those required for well-posedness. However, the first-order mixed-difference behaves as a product of one-direction errors~\citep{mimc2016} only when the solution map $(P_{\alpha_1},N_{\alpha_2}) \mapsto G_{\boldsymbol\alpha}$ has mixed regularity with respect to $P$ and $N$. This regularity ensures that $\abs{\E{{\hat{\Delta} G_{\boldsymbol\alpha}}}}$ approximately decays as a product of the geometric factors in each direction ($P$ and $N$), which is the mechanism behind sparse-grid/combination technique error bounds~\citep{pflaum1999error}. These rates can be numerically verified for a given problem. Section~\ref{sec:numerics} numerically investigates how the regularity of $G$ and $\kappa_1$ affects Assumption~\ref{ass:midlmc_bias} for the Kuramoto model~\eqref{eqn:kuramoto_model}. A full theoretical analysis of the conditions on $b,\sigma,\kappa_1,\kappa_2$ and $G$ ensuring Assumption~\ref{ass:midlmc_bias} is the subject of ongoing work.

The index set $\mathcal{I} \subset \mathbb{N}^2$ is selected to satisfy the bias constraint in~\eqref{eqn:midlmc_bias_constraint}. The statistical error constraint in~\eqref{eqn:midlmc_stat_constraint} is  expressed as the following constraint on the variance of the estimator using the asymptotic normality of the multi-index estimator~\citep{mimc2016}:

\begin{equation}
	\label{eqn:midlmc_var_constraint}
	\frac{\Var{\mathcal{A}_{\mathrm{MIMC}}(\mathcal{I})}}{\abs{\E{G(X(T))}}^2} \leq \left( \frac{\theta \tol_{\mathrm{r}}}{C_\nu} \right)^2,
\end{equation}
where $C_\nu$ denotes the $\left(1 - \frac{\nu}{2} \right)$-quantile of the standard normal distribution. Due to the independence of the DLMC estimators corresponding to each $\boldsymbol\alpha$ in \eqref{eqn:midlmc_estimator}, the variance of the multi-index DLMC estimator is written as follows:

\begin{equation}
	\label{midlmc_variance_v1}
	\Var{\mathcal{A}_{\mathrm{MIMC}}(\mathcal{I})} = \sum_{\boldsymbol\alpha \in \mathcal{I}} \frac{1}{M_{1,\boldsymbol\alpha}} \Var{\frac{1}{M_{2,\boldsymbol\alpha}} \sum_{m_2=1}^{M_{2,\boldsymbol\alpha}} {\hat{\Delta} G_{\boldsymbol\alpha}} \left( \omega_{1:P_{\alpha_1}}^{(\boldsymbol\alpha,\cdot)} \times {\tilde{\omega}^{(\boldsymbol\alpha,\cdot,m_2)}} \right)} \cdot
\end{equation}
Owing to the law of total variance,

\begin{align}
	\Var{\mathcal{A}_{\mathrm{MIMC}}(\mathcal{I})} &= \sum_{\boldsymbol\alpha \in \mathcal{I}} \frac{1}{M_{1,\boldsymbol\alpha}} \Bigg( \Bigg. \underbrace{ \Var{\E{{\hat{\Delta} G_{\boldsymbol\alpha}} \mid \omega_{1:P_{\alpha_1}}^{(\boldsymbol\alpha,\cdot)}}}}_{=V_{1,\boldsymbol\alpha}} \nonumber \\
	&\qquad + \frac{1}{M_{2,\boldsymbol\alpha}} \underbrace{\E{\Var{{\hat{\Delta} G_{\boldsymbol\alpha}} \mid \omega_{1:P_{\alpha_1}}^{(\boldsymbol\alpha,\cdot)}}}}_{=V_{2,\boldsymbol\alpha}} \Bigg. \Bigg) \nonumber \\
	\label{midlmc_variance_v2}
	&= \sum_{\boldsymbol\alpha \in \mathcal{I}} \left( \frac{V_{1,\boldsymbol\alpha}}{M_{1,\boldsymbol\alpha}} + \frac{V_{2,\boldsymbol\alpha}}{M_{1,\boldsymbol\alpha} M_{2,\boldsymbol\alpha}} \right) ,
\end{align}
where $V_{1,\boldsymbol\alpha}$ and $V_{2,\boldsymbol\alpha}$ are conditioned on the discretized empirical law $\mu^{P_{\alpha_1} \mid N_{\alpha_2}}$, described by the set of random variables $\omega_{1:P_{\alpha_1}}^{(\boldsymbol\alpha,\cdot)}$.

\begin{assump}{ (\textit{Multiplicative variance decay})}
\label{ass:midlmc_var}
	\hspace{1mm} For all $\boldsymbol\alpha \in \mathcal{I}$, for the constants $w_i > 0$ and $s_i >0$ for $i=1,2$, the variance terms $V_{1,\boldsymbol\alpha},V_{2,\boldsymbol\alpha}$ defined in \eqref{midlmc_variance_v2} satisfy
	\begin{align}
		\label{eqn:midlmc_v1_ass}
		V_{1,\boldsymbol\alpha} = \Var{\E{{\hat{\Delta} G_{\boldsymbol\alpha}} \mid \omega_{1:P_{\alpha_1}}^{(\boldsymbol\alpha,\cdot)}}} &\lesssim P_{\alpha_1}^{-w_1} N_{\alpha_2}^{-w_2}, \\
		\label{eqn:midlmc_v2_ass}
		V_{2,\boldsymbol\alpha} = \E{\Var{{\hat{\Delta} G_{\boldsymbol\alpha}} \mid \omega_{1:P_{\alpha_1}}^{(\boldsymbol\alpha,\cdot)}}} &\lesssim P_{\alpha_1}^{-s_1} N_{\alpha_2}^{-s_2} \cdot
	\end{align}
\end{assump}
Assumption~\ref{ass:midlmc_var} is motivated by the propagation of chaos~\citep{mean_field_limit} results for MV-SDEs and the strongly coupled antithetic estimator developed in~\eqref{eqn:antithetic_sampler}. The rates $w_2 = 2$ and $s_2 = 1$ are expected from the standard weak~\citep{talay1990expansion} and strong~\citep{sde_numerics} convergence results of the Euler--Maruyama scheme. Under the sufficient regularity of $b,\sigma,\kappa_1,\kappa_2$, and $G$ and antithetic coupling~\eqref{eqn:antithetic_coupling}, the weak~\citep{theory_mvsde} and strong~\citep{Szpruch:2021aa} convergence analysis of the particle system~\eqref{eqn:strong_approx_mvsde} indicate that $w_1 = 2$ and $s_1 = 2$. Under further mixed regularity, this coupling in the mixed-differences ensures that the particle approximation and time discretization errors combine multiplicatively in the variance bounds. With such a coupling scheme, the leading-order increments cancel in the mixed-difference so that $V_{1,\boldsymbol\alpha}$ and $V_{2,\boldsymbol\alpha}$ decay faster than the multilevel variances in the work by~\citep{my_mldlmcis}. These rates can be numerically checked for a given problem (see Section~\ref{sec:numerics} for the Kuramoto model). Section~\ref{sec:numerics} numerically examines how the regularity of $G$ and $\kappa_1$ influences  Assumption~\ref{ass:midlmc_var} for the Kuramoto model~\eqref{eqn:kuramoto_model}. A rigorous proof of these mixed convergence rates would require combining the propagation of chaos results with multi-index coupling arguments in~\eqref{eqn:antithetic_sampler} and is left for future work.

\subsection{Computational complexity}

For each multi-index $\boldsymbol\alpha \in \mathbb{N}^2$, $\E{{\hat{\Delta} G_{\boldsymbol\alpha}}}$ is estimated using six evaluations of $G_{\boldsymbol\alpha}$ for various combinations of discretization parameters~\eqref{eqn:antithetic_sampler}. The dominant computational contribution corresponds to the multi-index $\boldsymbol\alpha = (\alpha_1,\alpha_2)$. This work assumes a naive $\order{P}$ implementation to evaluate the empirical means in the drift and diffusion functions in~\eqref{eqn:strong_approx_mvsde} and~\eqref{eqn:dmvsde_sde_is}. The computational cost to generate one realization of $G_{\boldsymbol\alpha}$ comprises two parts: (i) generating  a realization of the empirical law $\mu^{P_{\alpha_1} \mid N_{\alpha_2}}$ using the Euler--Maruyama scheme with $N_{\alpha_2}$ time steps, requiring $\order{N_{\alpha_2} P_{\alpha_1}^2}$ operations, and (ii) generating a sample path of the decoupled MV-SDE~\eqref{eqn:dmvsde_sde_is} using the Euler--Maruyama scheme with $N_{\alpha_2}$ time steps, requiring $\order{N_{\alpha_2} P_{\alpha_1}}$ operations. The cost of computing the IS control $\zeta$ by numerically solving~\eqref{eqn:dmvsde_hjb_form3} is incurred only once, and the obtained IS control is employed uniformly across all multi-indices $\boldsymbol\alpha$; hence, it is treated as an off-line preprocessing step. Thus, the total computational cost of the multi-index DLMC estimator can be bound as follows:

\begin{equation}
	\label{eqn:midlmc_work}
	\mathcal{W}[\mathcal{A}_{\mathrm{MIMC}}(\mathcal{I})] \lesssim \sum_{\boldsymbol\alpha \in \mathcal{I}} M_{1,\boldsymbol\alpha} N_{\alpha_2} P_{\alpha_1}^{2} + M_{1,\boldsymbol\alpha} M_{2,\boldsymbol\alpha} N_{\alpha_2} P_{\alpha_1} \cdot
\end{equation}
Following the work by~\citep{mimc2016}, the optimal index set for every value of $L > 0$ under Assumptions~\ref{ass:midlmc_bias} and~\ref{ass:midlmc_var} can be written as follows:

\begin{equation}
\label{eqn:midlmc_optimal_indexset}
	\mathcal{I}(L) = \left\{ \boldsymbol\alpha \in \mathbb{N}^2 : (1-\bar{s}_1+2b_1)\alpha_1 + (1-\bar{s}_2+2b_2)\alpha_2 \leq L \right\},
\end{equation}
where $\bar{s}_1 = \min (w_1-1,s_1)$ and $\bar{s}_2 = \min (w_2,s_2)$. The following theorem presents the optimal computational complexity of the multi-index DLMC estimator with the above optimal index set~\eqref{eqn:midlmc_optimal_indexset}.

\begin{thm}{\textit{(Optimal multi-index DLMC complexity)}}
\label{thm:midlmc_complexity}
	\hspace{1mm} {Let $\theta \in (0,1)$ and $\nu \in (0,1)$ in~\eqref{eqn:midlmc_bias_constraint} and~\eqref{eqn:midlmc_var_constraint} be fixed.} Let Assumptions~\ref{ass:midlmc_bias} and~\ref{ass:midlmc_var} hold, along with the standard well-posedness assumptions for~\eqref{eqn:mvsde}. Let $b_1,b_2,w_1,w_2,s_1,s_2 > 0$ be the constants in Assumptions~\ref{ass:midlmc_bias} and~\ref{ass:midlmc_var}. {Let $\bar{L} \equiv \bar{L}(\tol_\mathrm{r})$ be chosen such that the bias constraint~\eqref{eqn:midlmc_bias_constraint} holds for the optimal multi-index set $\mathcal{I}(\bar{L})$ defined in~\eqref{eqn:midlmc_optimal_indexset}.}Let $\bar{s}_1 = \min (w_1-1,s_1)$ and $\bar{s}_2 = \min (w_2,s_2)$. If $2b_1 \geq \bar{s}_1-1$ and $2b_2 \geq \bar{s}_2-1$, then the total computational cost of the multi-index DLMC estimator~\eqref{eqn:midlmc_estimator} $\mathcal{W}[\mathcal{A}_{\mathrm{MIMC}}(\mathcal{I}(\bar{L}))]$, subject to the variance constraint~\eqref{eqn:midlmc_var_constraint}, satisfies
	\begin{equation}
	\label{eqn:midlmc_complexity}
		\limsup_{\tol_{\mathrm{r}} \downarrow 0} \frac{\mathcal{W}[\mathcal{A}_{\mathrm{MIMC}}(\mathcal{I}(\bar{L}))]}{\tol_{\mathrm{r}}^{-2-2 \max (0, \varsigma)} \left( \log \tol_{\mathrm{r}}^{-1} \right)^{\varrho}} \leq C_{\mathrm{work}} < \infty,
	\end{equation}
	where, 
	\begin{align*}
		\varsigma &= \max \left( \frac{1-\bar{s}_1}{2 b_1}, \frac{1-\bar{s}_2}{2 b_2} \right), \qquad \xi = \min \left( 2b_1 - \bar{s}_1, 2b_2 - \bar{s}_2 \right), \\
		\varrho &= \begin{cases}
		0, \quad \varsigma < 0 \\
		2 \aleph, \quad \varsigma = 0 \\
		1 + 2(\aleph - 1)(1 + \varsigma), \quad \varsigma>0 \text{ and } \xi = 0, \\
		2(\aleph-1)(1 + \varsigma), \quad \varsigma > 0 \text{ and } \xi > 0 \cdot
		\end{cases}, \qquad \aleph = \begin{cases}
		1, \quad \frac{1-\bar{s}_1}{2 b_1} \neq \frac{1-\bar{s}_2}{2 b_2}, \\
		2, \quad \frac{1-\bar{s}_1}{2 b_1} = \frac{1-\bar{s}_2}{2 b_2} \cdot
		\end{cases}
	\end{align*}
\end{thm}

\begin{proof}
	See~\ref{app:complexity}.
\end{proof}
For the Kuramoto model~\eqref{eqn:kuramoto_model} with the $C^\infty$-mollified indicator observable \\
$G(x) = \frac{1}{2} \left( 1 + \tanh (3(x-K)) \right)$, Section~\ref{sec:numerics} numerically demonstrates that $b_1=b_2=1$, $w_1=w_2=s_1=2$, and $s_2=1.5$. Using Theorem~\ref{thm:midlmc_complexity} with these values indicates that the complexity rate of the multi-index DLMC estimator is $\order{\tol_{\mathrm{r}}^{-2} \left( \log \tol_{\mathrm{r}}^{-1} \right)^2}$ compared with that of $\order{\tol_{\mathrm{r}}^{-3}}$ for the multilevel DLMC estimator~\citep{my_mldlmcis}, as confirmed in Figure~\ref{fig:adaptive_midlmcis}. This complexity matches the optimal rate predicted by multi-index theory for these convergence rates. Achieving a better asymptotic rate would likely require either higher-order methods or additional variance reduction beyond IS, which is left for future research.

\begin{rem}{\textit{(On isotropic directions)}}
	To demonstrate the superiority of the multi-index DLMC estimator over the multilevel DLMC estimator from the work by~\citep{my_mldlmcis}, this work examines the case in which $w_1=w_2=s_1=s_2=s>0$ and $b_1=b_2=b>0$. The expectations and variances of the mixed-differences converge at the same rate in the $P$ and $N$ directions. One such case is the Kuramoto setting~\eqref{eqn:kuramoto_model} with the smooth, nonrare observable $G(x) = \cos(x)$. Theorem~\ref{thm:midlmc_complexity}, as $\tol_{\mathrm{r}} \rightarrow 0$, yields the following result:
	\begin{equation}
	\label{eqn:midlmc_isotropic_complexity}
		\mathcal{W}[\mathcal{A}_{\mathrm{MIMC}}] = \begin{cases}
		\order{\tol_{\mathrm{r}}^{-2}}, \quad &s > 2 \\
		\order{\tol_{\mathrm{r}}^{-2} \left( \log \tol_{\mathrm{r}}^{-1} \right)^2}, \quad &s = 2 \\
		\order{\tol_{\mathrm{r}}^{-2-\frac{2-s}{b}}}, \quad &s < 2
		\end{cases} \cdot
	\end{equation}
	In contrast, the multilevel DLMC estimator~\citep{my_mldlmcis} has the following complexity rates:
	\begin{equation}
	\label{eqn:mldlmc_isotropic_complexity}
		\mathcal{W}[\mathcal{A}_{\mathrm{MLMC}}] = \begin{cases}
		\order{\tol_{\mathrm{r}}^{-2}}, \quad &s > 3 \\
		\order{\tol_{\mathrm{r}}^{-2} \left( \log \tol_{\mathrm{r}}^{-1} \right)^2}, \quad &s = 3 \\
		\order{\tol_{\mathrm{r}}^{-2-\frac{3-s}{b}}}, \quad &s < 3
		\end{cases} \cdot
	\end{equation}
	Comparing \eqref{eqn:midlmc_isotropic_complexity} and \eqref{eqn:mldlmc_isotropic_complexity}, the multi-index DLMC estimator has the same or better complexity rates than the multilevel DLMC estimator in all regimes of the variance convergence rate $s$.  However, the multi-index DLMC estimator requires mixed regularity (in the sense of Assumptions~\ref{ass:midlmc_bias} and \ref{ass:midlmc_var}) with respect to $P$ and $N$, whereas the multilevel DLMC estimator entails only ordinary regularity.
\end{rem}

\begin{rem}{(\textit{Computing probabilities)}}
\label{rem:probabilities}
The multi-index DLMC estimator~\eqref{eqn:midlmc_estimator} requires observables $G$ with sufficient mixed regularity in $P$ and $N$ (refer to Assumptions~\ref{ass:midlmc_bias} and~\ref{ass:midlmc_var}). An extension to rare-event probabilities of the form $\alpha_K = \prob{H(X(T)) > K}$, with $H:\mathbb{R}^d \rightarrow \mathbb{R}$ and $K$ in the tail of $H(X(T))$, corresponds to the discontinuous indicator observable $G(x) = \mathbbm{1}_{\{H(x)>K\}}$. The lack of regularity results in a poor expectation and variance decay of mixed-differences and {suboptimal complexity (Section~\ref{sec:numerics})}. A solution is to replace the indicator with {a mollifier $\tilde{\mathbbm{1}}^\varepsilon_{\{H(x)>K\}} \in C^\infty(\mathbb{R})$}, where $\varepsilon$ denotes the smoothing parameter. {If $H(X(T))$ admits a density that is sufficiently smooth in a neighbourhood of $K$, and if the mollification kernel is symmetric in the sense of~\citep{giles2015distribution}, then the absolute mollification bias satisfies}
	\begin{equation*}
		\abs{\E{\mathbbm{1}_{\{H(X(T))>K\}} - \tilde{\mathbbm{1}}^\varepsilon_{\{H(X(T))>K\}}}} = \order{\varepsilon^2} \cdot
	\end{equation*}
{Thus, setting $\varepsilon = \order{(\alpha_K \tol_\mathrm{r})^{0.5}}$ achieves a relative error tolerance $\tol_\mathrm{r}$~\citep{giles2015distribution,giles2023mlmc} and restores the complexity rates of Theorem~\ref{thm:midlmc_complexity}.} In practice, decreasing $\varepsilon$ requires an adaptive time-stepping method~\citep{szepessy2001adaptive} to maintain mixed-difference convergence rates because the mollified observable becomes increasingly steep as $\tol_\mathrm{r} \rightarrow 0$. The systematic development of such mollification-adaptivity strategies is the subject of ongoing work.
\end{rem}

\begin{rem}{(\textit{High-dimensional problems})}
\label{rem:high_dim}
	{For fixed dimension $d$ of the MV-SDE, Theorem~\ref{thm:midlmc_complexity} gives a complexity exponent in $\tol_\mathrm{r}$ that is independent of $d$. The constants and the per-sample cost, however, do depend on $d$, and the off-line cost of computing the IS control $\zeta$ by numerically solving the PDE~\eqref{eqn:dmvsde_hjb_form3} scales exponentially with $d$ using standard finite difference/element solvers. For higher-dimensional problems ($d \gg 1$), this step may become the dominant computational bottleneck. In such settings, the total computational cost should be viewed as $\mathcal{W}_\text{PDE} + \mathcal{W}[\mathcal{A}_\text{MIMC}]$, where $\mathcal{W}_\text{PDE}$ denotes the cost of numerically constructing the IS control $\zeta$. There is an inherent trade-off between these two contributions. Increasing the accuracy of the PDE solve requires finer discretisations and thus increases $\mathcal{W}_\text{PDE}$. However, the resulting control is closer to the exact IS control~\eqref{eqn:dmvsde_hjb_optimal_control}, leading to a smaller estimator variance and consequently less multi-index DLMC samples required to achieve a prescribed relative tolerance. Conversely, a less accurate PDE solve reduces the offline cost but typically increases the online sampling cost. Developing a framework that jointly optimizes the PDE accuracy in high dimensions and the multi-index DLMC sampling parameters, similar to the work in~\citep{BenAmar:2026aa}, is left for future research. 
To mitigate the curse of dimensionality, one may employ model reduction~\citep{Hartmann:2016aa}, variational formulation~\citep{Hartmann:2019aa}, sparse tensor techniques~\citep{Khoromskij:2015aa}, or machine learning approaches~\citep{Han:2018aa}. For the 1D Kuramoto model~\eqref{eqn:kuramoto_model} considered in this work, we observe that $\mathcal{W}_\text{PDE} \ll \mathcal{W}[\mathcal{A}_\text{MIMC}]$ over the range of prescribed relative error tolerances, and can hence treat the PDE solve as a preprocessing step.}
\end{rem}

\subsection{Adaptive algorithm}
\label{sec:adaptive_midlmcis}

Algorithm~\ref{alg:adaptive_midlmcis} presents the adaptive multi-index DLMC algorithm employed to construct the optimal index set $\mathcal{I}$~\eqref{eqn:midlmc_optimal_indexset} sequentially and determine the optimal number of samples $\{M_{1,\boldsymbol\alpha}, M_{2,\boldsymbol\alpha}\}_{\boldsymbol\alpha \in \mathcal{I}}$ to satisfy the relative error constraints in~\eqref{eqn:midlmc_bias_constraint} and \eqref{eqn:midlmc_stat_constraint}. This algorithm is adapted from the MIMC algorithm in the work by~\citep{mimc2016}. The algorithm employs pilot computations to estimate the convergence rates in Assumptions~\ref{ass:midlmc_bias} and~\ref{ass:midlmc_var} and allocates computational effort optimally across multi-indices via the optimal index set~\eqref{eqn:midlmc_optimal_indexset}. The adaptive strategy starts with a small index set (corresponding to $L=1$ in~\eqref{eqn:midlmc_optimal_indexset}), estimates the bias and variance contributions on those indices, and expands the index set until the bias criterion~\eqref{eqn:midlmc_bias_constraint} is met. At each expansion, this algorithm also computes the optimal number of samples per multi-index according to~\eqref{eqn:midlmc_actual_soln}. Therefore, the resulting computational cost of the adaptive algorithm scales as predicted by Theorem~\ref{thm:midlmc_complexity}. The algorithm relies on accurate, low-cost estimates of the relative bias and variances $\{V_{1,\boldsymbol\alpha}, V_{2,\boldsymbol\alpha}\}_{\boldsymbol\alpha \in \mathcal{I}}$. Small auxiliary simulations (Algorithms~\ref{alg:midlmcis_md}-\ref{alg:midlmc_extrapolation} in~\ref{app:algorithms}) estimate the bias/variances to guide adaptivity. The cost of these estimations is negligible compared with the total computational work as $\tol_\mathrm{r} \rightarrow 0$ because these estimations are performed with fixed sample sizes and interpolation techniques. Another critical feature of Algorithm~\ref{alg:adaptive_midlmcis} is that the IS control $\zeta$ is computed once off-line. A single realization of the empirical law $\mu^{\bar{P} \mid \bar{N}}$ is generated for large $\bar{P}$ and $\bar{N}$, after which the control PDE in~\eqref{eqn:dmvsde_hjb_form3}, given $\mu^{\bar{P} \mid \bar{N}}$, is solved numerically to obtain $\zeta$. \ref{app:algorithms} presents further details.

\begin{algorithm}
    \SetAlgoLined
    \textbf{User provided input}: $P_0,N_0,\tol_{\mathrm{r}},\nu,\theta,\{\bar{M}_1,\bar{M}_2\},\{\tilde{M}_1,\tilde{M}_2\}$; \\
    \textbf{Input from the numerical PDE solver}: IS control $\zeta(\cdot,\cdot)$; \\
    \textbf{Input from pilots}: $\{b_1,b_2\},\{w_1,w_2\},\{s_1,s_2\}$; \\
    Estimate $\bar{G} = \E{{\hat{\Delta} G_{(0,0)}}}$ using \textbf{Algorithm~\ref{alg:midlmcis_md}} with $\bar{M}_1,\bar{M}_2,\zeta(\cdot,\cdot)$; \\
    Estimate and store $\{V_{1,\boldsymbol\alpha},V_{2,\boldsymbol\alpha}\}$ for $\boldsymbol\alpha \in [0,1,2] \cross [0,1,2]$ using \textbf{Algorithm~\ref{alg:adaptive_variance}} with $\tilde{M}_1,\tilde{M}_2,\zeta(\cdot,\cdot)$; \\
    Set $L = 1$ {and initialize $\hat{\epsilon}_b = +\infty$}; \\
    {
    \Repeat{$\hat{\epsilon}_b > (1-\theta) \tol_{\mathrm{r}}$}
    {
    Generate index set $\mathcal{I}(L)$ from \eqref{eqn:midlmc_optimal_indexset}; \\
    Estimate and store $\{V_{1,\boldsymbol{\alpha}},V_{2,\boldsymbol{\alpha}}\}_{\boldsymbol{\alpha} \in \mathcal{I}(L)}$ using \textbf{Algorithm~\ref{alg:midlmc_extrapolation}}; \\
    Compute the optimal sample sizes $\{M_{1,\boldsymbol{\alpha}},M_{2,\boldsymbol{\alpha}}\}_{\boldsymbol{\alpha} \in \mathcal{I}(L)}$ using \eqref{eqn:midlmc_actual_soln};\\
    Reevaluate $\bar{G} =  \sum_{\boldsymbol{\alpha} \in \mathcal{I}(L)} \E{{\hat{\Delta} G_{\boldsymbol\alpha}}}$ as per \eqref{eqn:midlmc_estimator} with $\{M_{1,\boldsymbol{\alpha}},M_{2,\boldsymbol{\alpha}}\},\zeta(\cdot,\cdot)$ using \textbf{Algorithm~\ref{alg:midlmcis_md}} for each $\boldsymbol{\alpha}$; \\
    $\hat{\epsilon}_b = \frac{1}{\abs{\bar{G}}} \sum_{\boldsymbol\alpha \in \partial \mathcal{I}(L)} \abs{\E{{\hat{\Delta} G_{\boldsymbol\alpha}}}}$; \\
    $L \leftarrow L+1$; \\
    }
    }
    $\mathcal{A}_{\mathrm{MIMC}} = \bar{G}$. \\
    \caption{Adaptive multi-index double loop Monte Carlo estimator with importance sampling}
    \label{alg:adaptive_midlmcis}
\end{algorithm}

\section{Numerical results}
\label{sec:numerics}

\begin{figure}[h!]
	\centering
	\includegraphics[width=.7\textwidth]{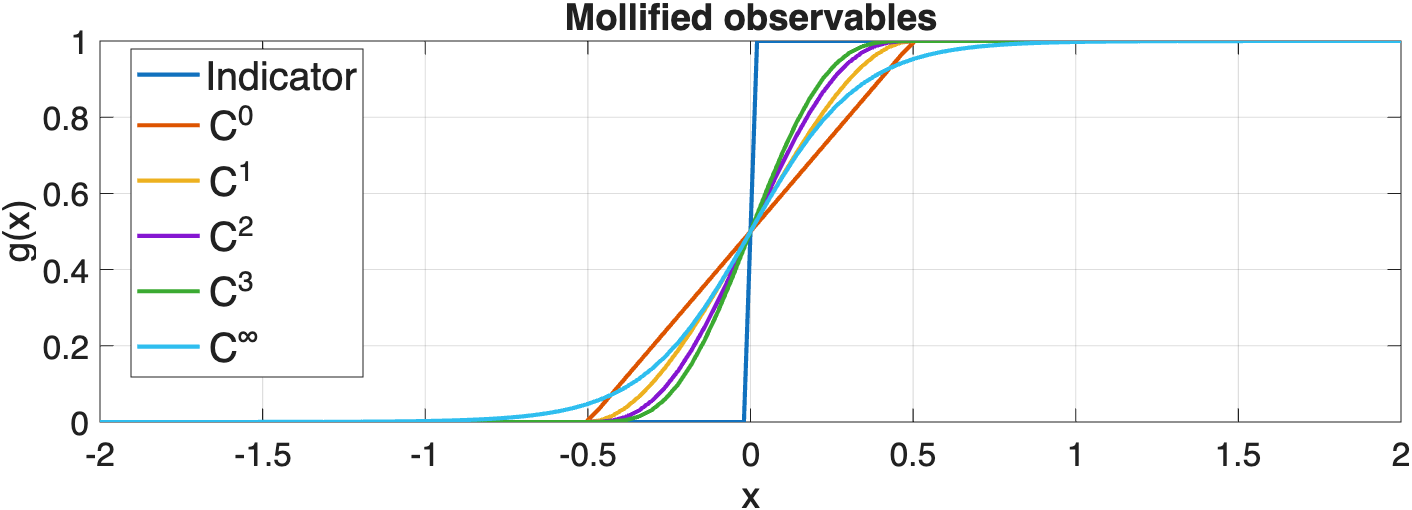}
	\caption{Rare-event observables $G$ plotted for $\epsilon = 0.5$ and $K=0$, including the indicator function and corresponding mollifiers of varying regularity.}
	\label{fig:mollifiers}
\end{figure}

This section presents numerical evidence supporting the assumptions and complexity rates derived in Section~\ref{sec:midlmc} and numerically verifies the variance reduction achieved using the IS scheme. All experiments\footnote{Numerical results were produced using MATLAB (v2022b) on an Intel Xeon Platinum 8360Y CPU (2.40~GHz).} were performed on the Kuramoto model~\eqref{eqn:kuramoto_model} with the following parameters: $\sigma = 0.4$, $T = 1$, $x_0 \sim \mathcal{N}(0,0.2)$, and $\xi \sim \mathcal{U}(-0.2,0.2)$. This work sets the tolerance-splitting parameter to $\theta=0.5$ and the prescribed confidence level to $\nu=0.05$ and defines the following:

\begin{equation}
\label{eqn:midlmc_hierarchy}
	P_{\alpha_1} = 5 \times 2^{\alpha_1}, \quad N_{\alpha_2} = 4 \times 2^{\alpha_2} .
\end{equation}
The multi-index DLMC method was assessed on the indicator observable $G(x) = \mathbbm{1}_{\{x > K\}}$ and on several mollified variants with varying regularities, where $K \in \mathbb{R}$ denotes a threshold parameter determining the event rarity. Each mollifier is specified by a  smoothing parameter $\varepsilon$ and is defined as follows:

\begin{equation}
\label{eqn:poly_mollifiers}
	G(x) = \begin{cases}
	0, \quad &x \leq K-\varepsilon\\
	S \left( \frac{x-K+\varepsilon}{2 \varepsilon} \right), \quad &K-\varepsilon < x \leq K+\varepsilon\\
	1, \quad &x > K+\varepsilon
	\end{cases},
\end{equation}
where $S$ denotes a polynomial selected to enforce the desired regularity. Table~\ref{tab:mollifiers} lists the explicit forms. {This work also considers the mollifier $G(x) = \frac{1}{2} \left(1 + \tanh \left( \frac{x-K}{\varepsilon} \right) \right)$ where $G \in C^\infty(\mathbb{R})$.} Figure~\ref{fig:mollifiers} depicts all observables in this section. For each choice of observable $G$, the IS control $\zeta$ from~\eqref{eqn:dmvsde_hjb_optimal_control} is computed by simulating a particle system~\eqref{eqn:strong_approx_mvsde} off-line with $\bar{P} = 1000$ particles and $\bar{N} = 100$ time steps. These values were selected for the control computation so that the suboptimality of the control is negligible. Numerical tests found that further increasing this $\bar{P},\bar{N}$ did not noticeably alter the results, indicating that the choice is sufficient. The PDE~\eqref{eqn:dmvsde_hjb_form3} is numerically solved using finite differences with linear interpolation.

\begin{table}[h!]
\centering
	\begin{tabular}{|c|c|}
	\hline 
	Mollifier regularity & $S(x)$ \\ 
	\hline
	$C^0(\mathbb{R})$ & $x$ \\ 
	$C^1(\mathbb{R})$ & $-2x^3 + 3x^2$ \\  
	$C^2(\mathbb{R})$ & $6x^5 - 15x^4 + 10x^3$ \\
	$C^3(\mathbb{R})$ & $-20x^7 + 70x^6 - 84x^5 + 35x^4$\\
	\hline 
	\end{tabular}
	\caption{List of polynomial functions applied in~\eqref{eqn:poly_mollifiers} to obtain mollifiers of varying regularity.}
	\label{tab:mollifiers}
\end{table}

\subsection{Variance reduction}
\label{sec:numerics_variance}

\begin{figure}
    \centering
    \begin{subfigure}[b]{0.45\textwidth}
        \centering
        \includegraphics[width=\textwidth]{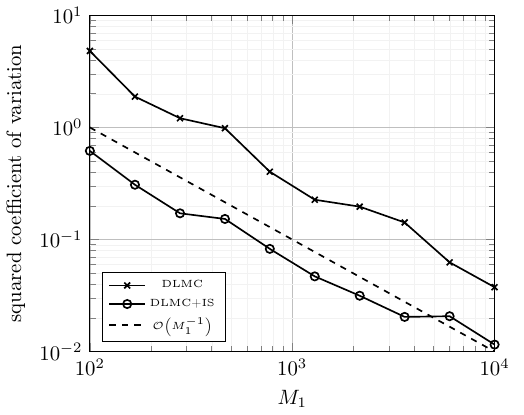}
        \caption{Squared coefficient of variation of the DLMC estimator for $\E{{\hat{\Delta} G_{(2,2)}}}$ with and without IS versus $M_1$ for fixed $M_2=100$.}
        \label{fig:is_cov}
    \end{subfigure}
    \hfill
    \begin{subfigure}[b]{0.45\textwidth}
        \centering
        \includegraphics[width=\textwidth]{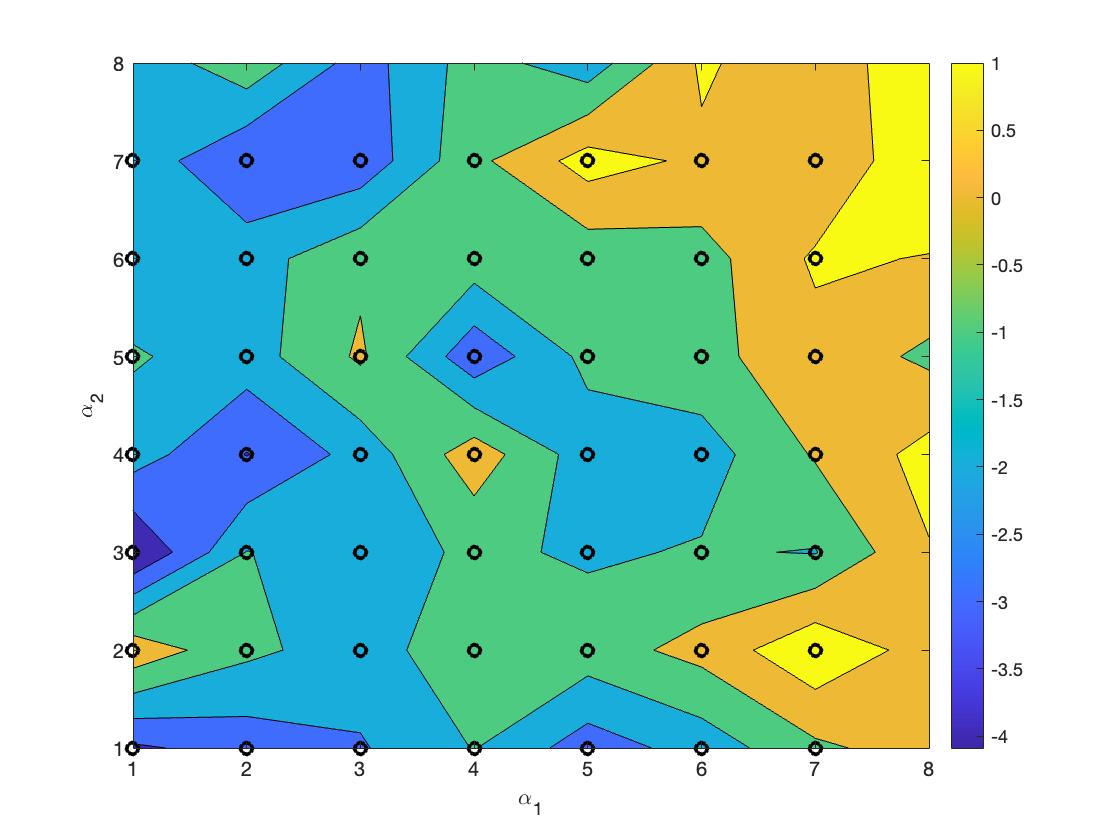}
        \caption{Contour plot (log scale) of $\mathcal{R}$ for a set of $\boldsymbol\alpha \in \mathbb{N}^2$ for the DLMC estimator of $\E{{\hat{\Delta} G_{\boldsymbol\alpha}}}$, presenting significant variance reduction, especially in coarser multi-indices.}
        \label{fig:is_contour_md}
    \end{subfigure}
    \caption{Variance reduction due to the IS scheme introduced in Section~\ref{sec:dlmcis} for the Kuramoto model~\eqref{eqn:kuramoto_model} with the $C^\infty$-mollifier observable with $\varepsilon = \frac{1}{3}$ and $K = 3.5$.}
    \label{fig:is_var_red}
\end{figure}

First, this work verifies variance reduction for DLMC estimators of $\E{{\hat{\Delta} G_{\boldsymbol\alpha}}}$ for the $C^\infty$-mollifier with $\varepsilon = \frac{1}{3}$ and $K=3.5$, for which $\E{G(X(T))} \approx 2.04 \cross 10^{-5}$. Figure~\ref{fig:is_cov} reports the squared coefficient of variation ($\frac{\Var{{\hat{\Delta} G_{\boldsymbol\alpha}}}}{\abs{\E{{\hat{\Delta} G_{\boldsymbol\alpha}}}}^2}$) as a function of $M_1$ with fixed $M_2$ for a given $\boldsymbol\alpha$, comparing DLMC with and without IS. For $\boldsymbol\alpha = (2,2)$, IS reduces the {coefficient of variation} by roughly an order of magnitude, while both estimators retain the expected $\order{M_1^{-1}}$ decay. Then, this work  examines the effect of the IS scheme across all multi-indices $\boldsymbol\alpha \in \mathbb{N}^2$, via the following ratio:
\begin{equation}
	\mathcal{R} = \frac{\Var{\mathcal{A}_{\mathrm{IS}}}}{\Var{\mathcal{A}_{\mathrm{MC}}}},
\end{equation}
where $\mathcal{A}_{\mathrm{IS}}$ and $\mathcal{A}_{\mathrm{MC}}$ denote the {DLMC estimators of $\E{{\hat{\Delta} G_{\boldsymbol\alpha}}}$ with and without IS}, respectively. Figure~\ref{fig:is_contour_md} presents a contour plot of $\mathcal{R}$ across various $\boldsymbol\alpha \in \mathbb{N}^2$. Up to $\order{10^3}$ variance reduction was observed on the coarse multi-indices, demonstrating the effectiveness of the IS control $\zeta$ when combined with the multi-index DLMC estimator in addressing the explosive growth of the coefficient of variation. This variance reduction directly enables the success of the multi-index method for rare-event quantities.

\subsection{Verifying Assumptions~\ref{ass:midlmc_bias} and~\ref{ass:midlmc_var}}
\label{sec:numerics_assumptions}

\begin{figure}
    \centering
    \begin{subfigure}[b]{0.6\textwidth}
        \centering
        \includegraphics[width=\textwidth]{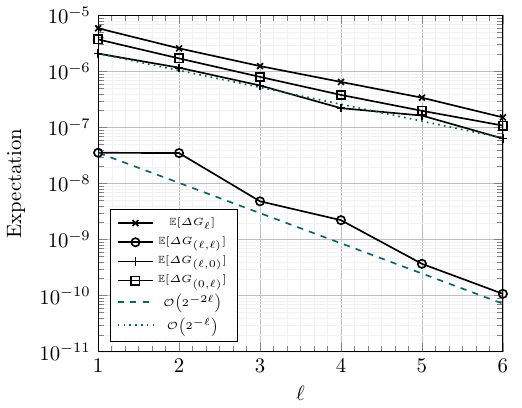}
        \caption{Verification of \eqref{eqn:midlmc_bias_ass} from Assumption~\ref{ass:midlmc_bias} using Algorithm~\ref{alg:midlmcis_md}, where $M_1=10^3$ and $M_2=10^3$.}
        \label{fig:midlmcis_rate_bias}
    \end{subfigure}
    \hfill
    \begin{subfigure}[b]{0.45\textwidth}
        \centering
        \includegraphics[width=\textwidth]{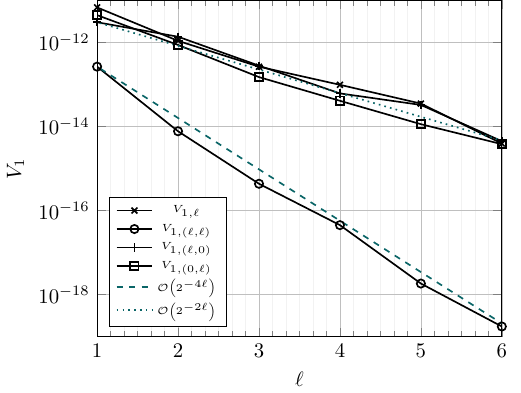}
        \caption{Verification of \eqref{eqn:midlmc_v1_ass} from Assumption~\ref{ass:midlmc_var} using Algorithm~\ref{alg:adaptive_variance}, where $M_1=10^2$ and $M_2=10^4$.}
        \label{fig:midlmcis_rate_v1a}
    \end{subfigure}
    \hfill
    \begin{subfigure}[b]{0.45\textwidth}
        \centering
        \includegraphics[width=\textwidth]{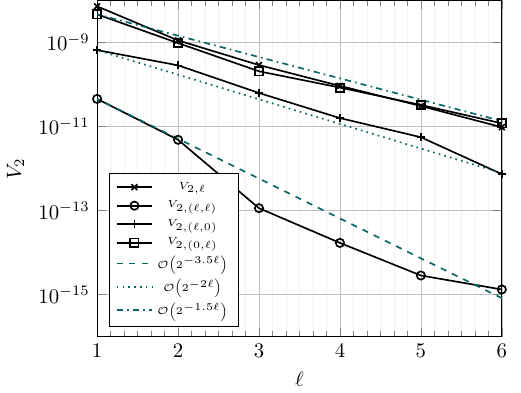}
        \caption{{Verification of \eqref{eqn:midlmc_v2_ass} from Assumption~\ref{ass:midlmc_var} using Algorithm~\ref{alg:adaptive_variance},} where $M_1=10^2$ and $M_2=10^4$.}
        \label{fig:midlmcis_rate_v2a}
    \end{subfigure}
    \caption{Kuramoto example, rate verification ($C^\infty$-mollifier observable with $K = 3.5$ and $\varepsilon = \frac{1}{3}$): Sample means and variances of mixed-differences in the multi-index DLMC estimator~\eqref{eqn:midlmc_estimator}.}
    \label{fig:midlmcis_rate}
\end{figure}

\begin{figure}
    \centering
    \begin{subfigure}[b]{0.6\textwidth}
        \centering
        \includegraphics[width=\textwidth]{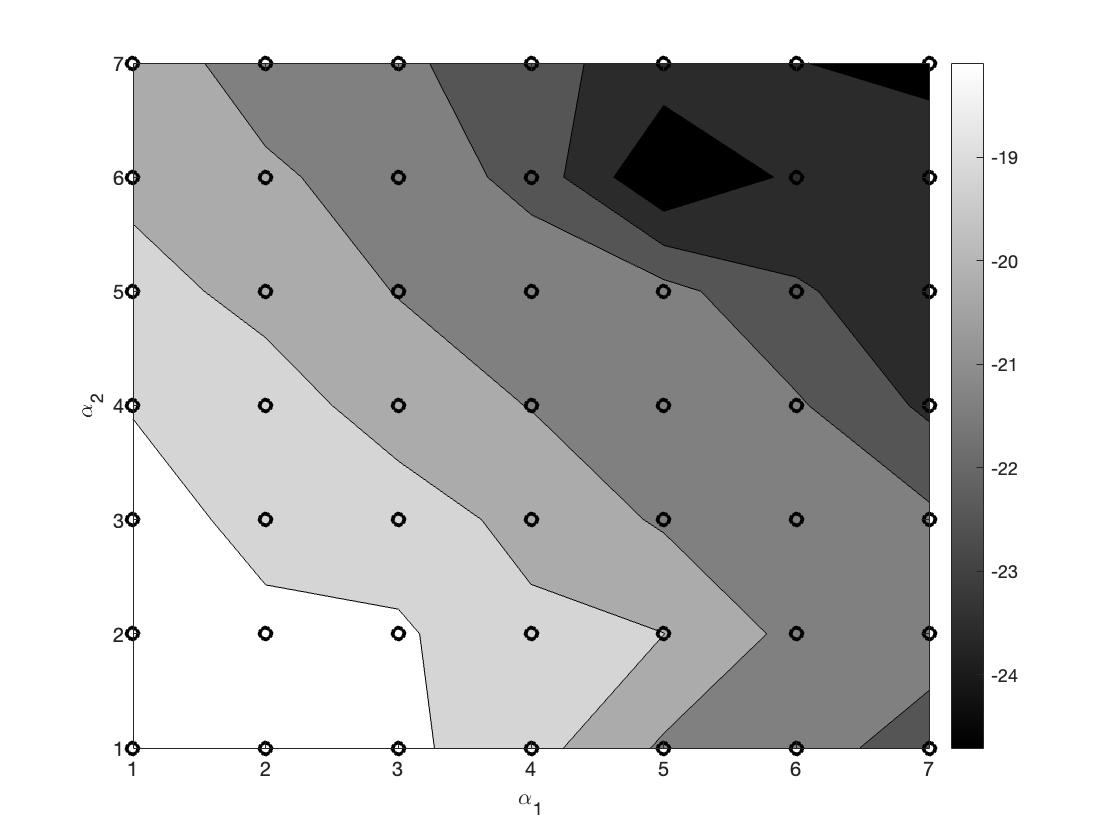}
        \caption{Contour plots of the sample mean $\abs{\E{{\hat{\Delta} G_{\boldsymbol\alpha}}}}$ of the mixed-differences used in the multi-index DLMC estimator~\eqref{eqn:midlmc_estimator} obtained using Algorithm~\ref{alg:midlmcis_md}, where $M_1=10^3$ and $M_2=10^2$.}
        \label{fig:md_exp_contour}
    \end{subfigure}
    \hfill
    \begin{subfigure}[b]{0.45\textwidth}
        \centering
        \includegraphics[width=\textwidth]{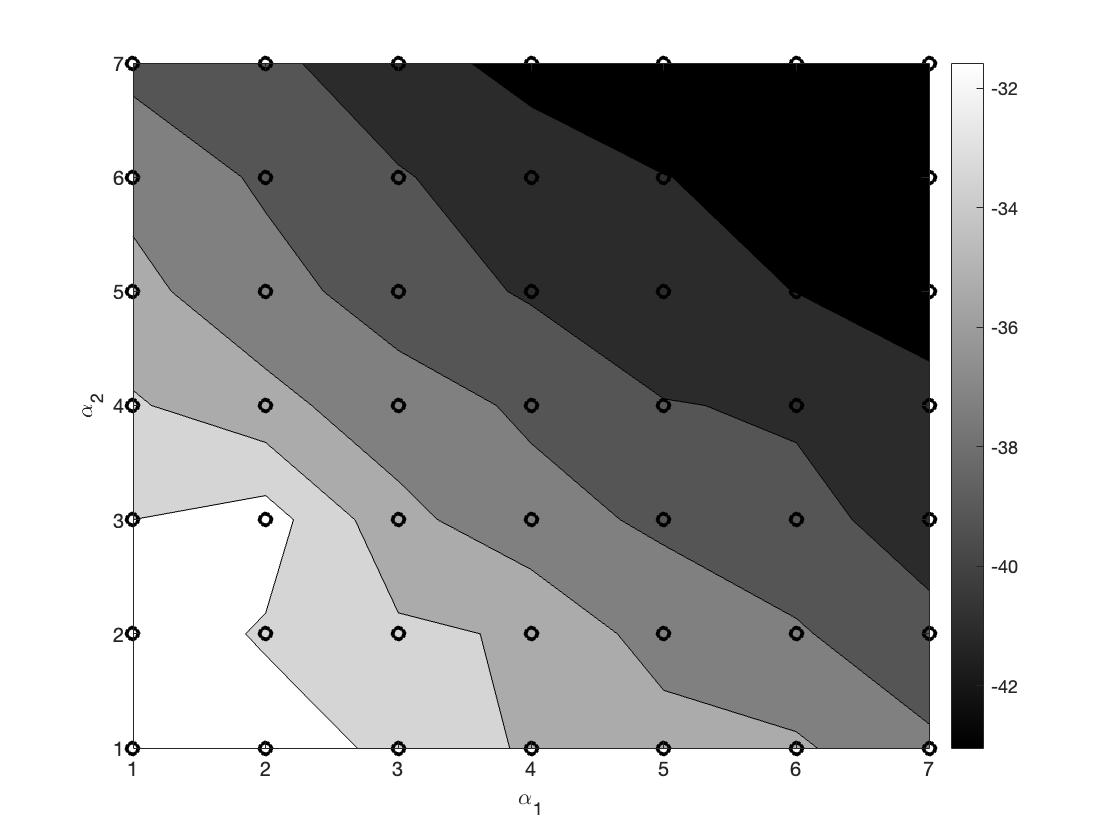}
        \caption{Contour plots of sample $V_{1,\boldsymbol\alpha}$ of the mixed-differences used in the multi-index DLMC estimator~\eqref{eqn:midlmc_estimator} obtained using Algorithm~\ref{alg:adaptive_variance}, where $M_1=10^2$ and $M_2=10^3$.}
        \label{fig:md_v1_contour}
    \end{subfigure}
    \hfill
    \begin{subfigure}[b]{0.45\textwidth}
        \centering
        \includegraphics[width=\textwidth]{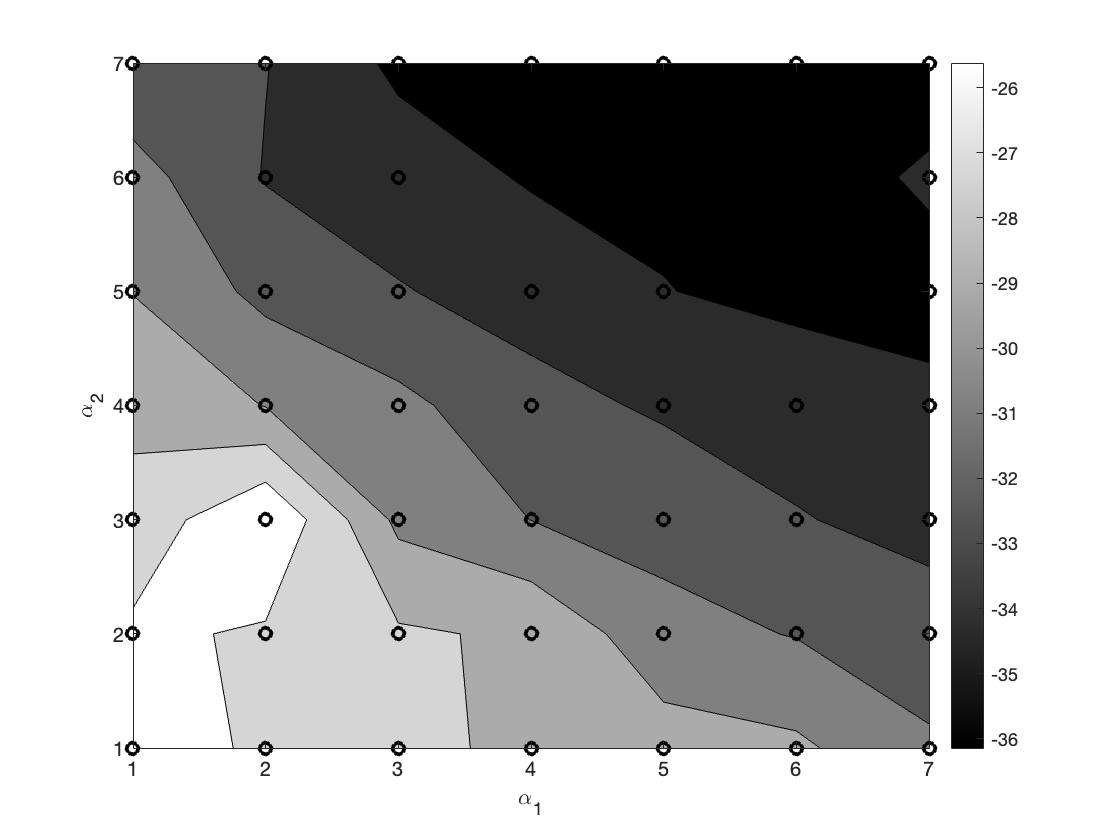}
        \caption{Contour plots of sample $V_{2,\boldsymbol\alpha}$ of the mixed-differences used in the multi-index DLMC estimator~\eqref{eqn:midlmc_estimator} obtained {using Algorithm~\ref{alg:adaptive_variance}}, where $M_1=10^2$ and $M_2=10^3$.}
        \label{fig:md_v2_contour}
    \end{subfigure}
    \caption{Kuramoto example, numerical rate verification ($C^\infty$-mollifier observable with $K = 3.5$ and $\varepsilon = \frac{1}{3}$): Contour plots of the sample mean and variances in the multi-index DLMC estimator~\eqref{eqn:midlmc_estimator}. All plots numerically and asymptotically verify Assumptions~\ref{ass:midlmc_bias} and \ref{ass:midlmc_var}.}
    \label{fig:midlmcis_contour}
\end{figure}

Figure~\ref{fig:midlmcis_rate} presents numerical evidence for  Assumptions~\ref{ass:midlmc_bias} and \ref{ass:midlmc_var} and estimates  the convergence rates $\{b_1,b_2\}$, $\{w_1,w_2\}$, and $\{s_1,s_2\}$ for the Kuramoto example with the same $C^\infty$-mollifier as above. Figure~\ref{fig:midlmcis_rate_bias} presents the decay of $\abs{\E{{\hat{\Delta} G_{\boldsymbol\alpha}}}}$, yielding rates $b_1 = 1$ and $b_2 = 1$. The figure also confirms that the mixed-difference expectation along the diagonal $(\ell,\ell)$ converges at a rate of $b_1+b_2=2$, demonstrating that it decays as the product of the single-direction errors. Figures ~\ref{fig:midlmcis_rate_v1a} and \ref{fig:midlmcis_rate_v2a} display the convergence of $V_{1, \boldsymbol\alpha}$ and $V_{2, \boldsymbol\alpha}$, respectively, yielding the rates $w_1 = 2$, $w_2 = 2$, $s_1 = 2$, and $s_2 = 1.5$. These results confirm the multiplicative decay of the mixed-difference variances. For comparison, we also plot the level differences $\E{\Delta G_\ell} = \E{G_{\ell} - G_{\ell-1}}$ from the multilevel DLMC estimator~\citep{my_mldlmcis}. The mixed-differences display higher convergence rates than the level differences of the multilevel estimator. Figure~\ref{fig:midlmcis_contour} reveals that Assumptions~\ref{ass:midlmc_bias} and ~\ref{ass:midlmc_var} are satisfied for sufficiently fine discretizations.

For further assessment of the validity of Assumptions~\ref{ass:midlmc_bias} and~\ref{ass:midlmc_var}, the same plots were repeated for observables $G$ of varying regularity, ranging from the discontinuous indicator function to the $C^3$-mollifier. Table~\ref{tab:regularity_rates} reports the numerical convergence rates in the $P$-direction $(\ell,0)$, the $N$-direction $(0,\ell)$, and the diagonal $(\ell,\ell)$-direction. \ref{app:numerics_addendum} presents the corresponding plots. For observables that are {in $C^k(\mathbb{R})$ for $k \geq 2$}, first-order bias decay occurs in both $P$ and $N$ directions, and second-order $V_{1,\boldsymbol\alpha}$ decay occurs in $P$ and $N$. {For the Kuramoto example and the tested discretisation range, $C^2$ regularity appears sufficient to recover the  multiplicative decay of the mixed-difference expectations and variances needed for Assumptions~\ref{ass:midlmc_bias} and~\ref{ass:midlmc_var}.} The multiplicative bias and variance convergence rates degrade with the $C^0$-mollifier, whereas slower convergence occurs in individual directions $P$ and $N$ for the indicator observable. Moreover, the numerical experiments demonstrate that the mixed-difference variances do not decay multiplicatively when replacing the $C^\infty$, bounded sinusoidal kernel in the Kuramoto model~\eqref{eqn:kuramoto_model} with a $C^0$, bounded triangular wave kernel, indicating that sufficient drift regularity is also necessary to satisfy Assumptions~\ref{ass:midlmc_bias} and~\ref{ass:midlmc_var}. This result confirms that sufficient smoothness in the MV-SDE~\eqref{eqn:mvsde} and observable $G$ is necessary for the feasibility of the multi-index DLMC estimator. A theoretical explanation of this behaviour is the subject of ongoing work. 

According to Theorem~\ref{thm:midlmc_complexity}, the mixed-difference convergence rates obtained for observables with $C^2$-regularity (and higher) ensure that the optimal complexity of the multi-index DLMC estimator is $\order{\tol_{\mathrm{r}}^{-2} \left( \log \tol_{\mathrm{r}}^{-1} \right)^2}$, which is one order better than the optimal work complexity of $\order{\tol_{\mathrm{r}}^{-3}}$ for the multilevel DLMC estimator~\citep{my_mldlmcis}.  The empirically observed rates $(b_i,w_i,s_i)_{i=1,2}$ validate Assumptions~\ref{ass:midlmc_bias} and~\ref{ass:midlmc_var} for sufficiently smooth cases, supporting the use of the complexity predictions of Theorem~\ref{thm:midlmc_complexity}.

\begin{table}[h!]
	\centering
	\begin{tabular}{||c|c|c|c|c|c|c|c|c|c||}
	\hline 
	$G(x)$ & \multicolumn{3}{c|}{$\abs{\E{{\hat{\Delta} G_{\boldsymbol\alpha}}}}$} &  \multicolumn{3}{c|}{$V_{1,\boldsymbol\alpha}$} &  \multicolumn{3}{c|}{$V_{2,\boldsymbol\alpha}$} \\ 
	\hline 
	 & $(\ell,\ell)$ & $(\ell,0)$ & $(0,\ell)$ & $(\ell,\ell)$ & $(\ell,0)$ & $(0,\ell)$ & $(\ell,\ell)$ & $(\ell,0)$ & $(0,\ell)$  \\ 
	\hline 
	$\mathbbm{1}_{\{x >K\}}$ & 0.7 & 1.1 & 0.3 & 1.2 & 1.4 & 1 & 1.1 & 0.6 & 0.6 \\ 
	\hline 
	$C^0$-mollifier & 1.8 & 1 & 1.6 & 3.1 & 2.1 & 1.7 & 2 & 1.8 & 0.8 \\ 
	\hline 
	$C^1$-mollifier & 1.9 & 1 & 1 & 3.6 & 2.1 & 1.8 & 2.5 & 2 & 1 \\ 
	\hline 
	$C^2$-mollifier & 2 & 1 & 0.8 & 3.7 & 2.1 & 1.8 & 2.6 & 2 & 1 \\ 
	\hline 
	$C^3$-mollifier & 2.2 & 1 & 0.8 & 3.8 & 2.1 & 1.8 & 2.6 & 2 & 1 \\ 
	\hline 
	$C^\infty$-mollifier & 2 & 1 & 1 & 4 & 2 & 2 & 3.5 & 2 & 1.5 \\ 
	\hline 
	\end{tabular}
	\caption{Comparison of expectation and variance convergence rates of mixed-differences for the Kuramoto model~\eqref{eqn:kuramoto_model} with mollifiers of varying regularity of the indicator function $G(x) = \mathbbm{1}_{\{x>K\}}$ for $K = 2$ and $\varepsilon = 0.5$.}
	\label{tab:regularity_rates} 
\end{table}

\subsection{Computational complexity: Multi-index vs. multilevel}
\label{sec:numerics_complexity}
\begin{figure}
    \centering
    \begin{subfigure}[b]{0.45\textwidth}
        \centering
        \includegraphics[width=\textwidth]{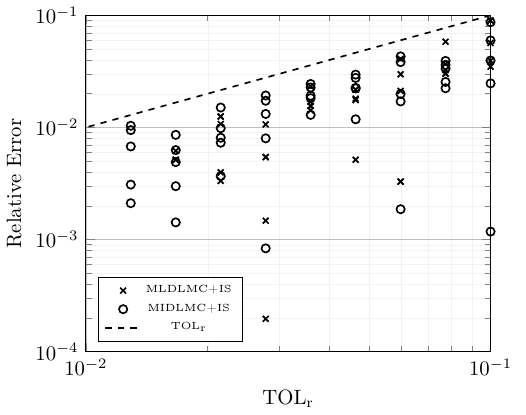}
        \caption{Estimated relative error for multilevel and multi-index DLMC estimators with IS.}
        \label{fig:global_err}
    \end{subfigure}
    \hfill
    \begin{subfigure}[b]{0.45\textwidth}
        \centering
        \includegraphics[width=\textwidth]{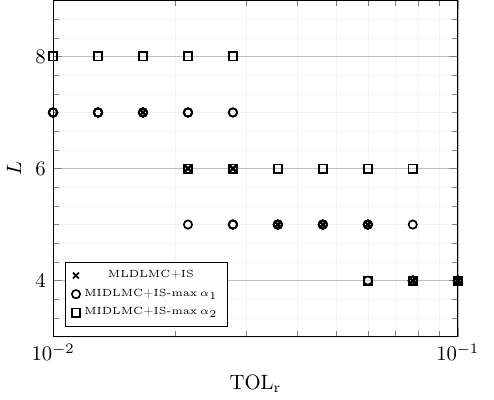}
        \caption{Maximum discretization level of the number of time steps and particles for the multi-index and multilevel DLMC estimators for different tolerances.}
        \label{fig:level}
    \end{subfigure}
    \hfill
    \begin{subfigure}[b]{0.45\textwidth}
        \centering
        \includegraphics[width=\textwidth]{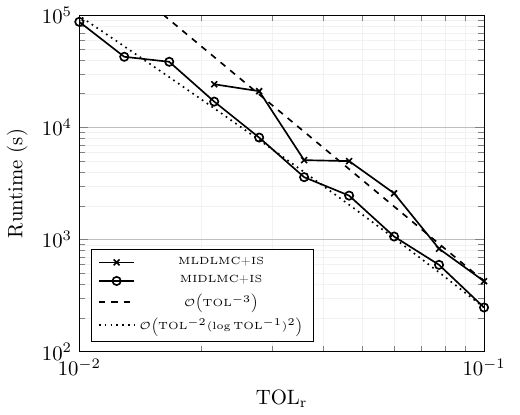}
        \caption{Average computational runtime for different values of $\tol_{\mathrm{r}}$.}
        \label{fig:runtime}
    \end{subfigure}
    \hfill
    \begin{subfigure}[b]{0.45\textwidth}
        \centering
        \includegraphics[width=\textwidth]{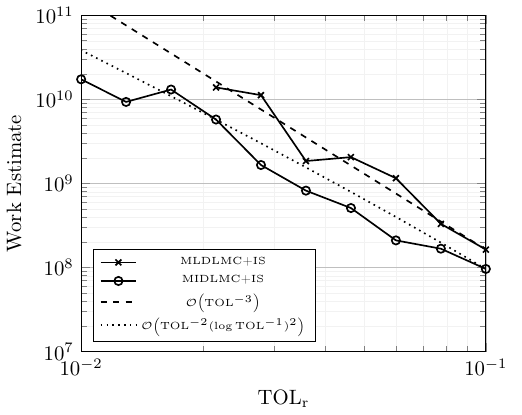}
        \caption{Average computational cost estimate~\eqref{eqn:midlmc_comp_cost} for different values of $\tol_{\mathrm{r}}$.}
        \label{fig:work_est}
    \end{subfigure}
    \caption{Kuramoto example, adaptive DLMC algorithms with IS ($C^\infty$-mollifier observable with $K = 3.5$ and $\varepsilon = \frac{1}{3}$): Comparing multilevel and multi-index estimators. Here, MLDLMC+IS denotes multilevel DLMC with IS and MIDLMC+IS denotes multi-index DLMC with IS.}
    \label{fig:adaptive_midlmcis}
\end{figure}

Figure~\ref{fig:adaptive_midlmcis} illustrates the results of running the adaptive multi-index DLMC Algorithm~\ref{alg:adaptive_midlmcis} for the $C^\infty$-mollifier with $K=3.5$ and $\varepsilon=\frac{1}{3}$, comparing them with the results from the adaptive multilevel algorithm in the work by~\citep{my_mldlmcis}. $K = 3.5$ corresponds to $\E{G(X(T))} \approx 2.04 \times 10^{-5}$, a sufficiently rare event. This work sets $\{\bar{M}_1, \bar{M}_2\} = \{10^3, 10^2\}$ and $\{\tilde{M}_1, \tilde{M}_2\} = \{25, 100\}$ (Algorithm~\ref{alg:adaptive_midlmcis}). Algorithm~\ref{alg:adaptive_midlmcis} and the corresponding multilevel algorithm were independently executed five times, and the combined results were plotted. Both algorithms apply the same IS control; the only difference is the hierarchy of discretizations (multi-index vs multilevel). 

Figure~\ref{fig:global_err} displays the exact relative error achieved using the multi-index and multilevel DLMC estimators for various relative tolerances $\tol_{\mathrm{r}}$, estimated using a reference multi-index DLMC approximation with $\tol_{\mathrm{r}} = 1\%$. Each marker represents a separate iteration of the corresponding adaptive algorithm. Figure~\ref{fig:global_err} illustrates the effectiveness of IS in estimating a rare-event quantity within a relative error tolerance, which is infeasible with the standard MC method. The result confirms that the multi-index and multilevel estimators with IS satisfy $\tol_{\mathrm{r}}$, in the sense of \eqref{eqn:objective}. 

Figure~\ref{fig:level} displays the maximum discretization level for the number of particles and time steps to ensure that the multilevel and multi-index estimators satisfy the relative bias constraint~\eqref{eqn:midlmc_bias_constraint}. For $\tol_\mathrm{r} = 1.67\%$, the multilevel DLMC estimator required $P = 640$ particles and $N = 512$ time steps, whereas the multi-index DLMC estimator increased to $P = 640$ particles and $N = 1024$ time steps with a lower net computational work. This result highlights the efficiency of the MIMC method over the multilevel MC method in distributing computational work across multi-indices, while satisfying the same relative accuracy.

Figure~\ref{fig:runtime} presents the average computational runtime for Algorithm~\ref{alg:adaptive_midlmcis} for various relative error tolerances, confirming the asymptotic complexity rates derived in Theorem~\ref{thm:midlmc_complexity}. The multilevel DLMC runtime scales with a power of roughly $-3$, whereas the multi-index DLMC scales at nearly $-2$. Moreover, the multi-index estimator not only asymptotically outperforms the multilevel method, but the absolute runtime is significantly lower for all tested tolerances due to the considerable IS  variance reduction (smaller constant in complexity).

Figure~\ref{fig:work_est} presents the average computational cost for both estimators for different values of $\tol_{\mathrm{r}}$, computed using the following computational cost model~\eqref{eqn:midlmc_work}:

\begin{equation}
	\label{eqn:midlmc_comp_cost}
	\text{Computational Cost}[\mathcal{A}_{\mathrm{MIMC}}(\mathcal{I}(L))] \approx \sum_{\boldsymbol\alpha \in \mathcal{I}(L)} \left( M_{1,\boldsymbol\alpha} N_{\alpha_2} P_{\alpha_1}^{2} + M_{1,\boldsymbol\alpha} M_{2,\boldsymbol\alpha} N_{\alpha_2} P_{\alpha_1} \right) \cdot
\end{equation}
For $\tol_\mathrm{r} = 2.78\%$, the multi-index DLMC estimator required $\sim 2 \times 10^{9}$ operations on average whereas the multilevel estimator required $\sim 10^{10}$ operations. For perspective, a direct MC simulation (without IS or multilevel/multi-index MC) requires $\sim 10^{15}$ operations to achieve the same relative accuracy given the event  rarity. 

In conclusion, this numerical study confirms the theoretical convergence rates and complexity derived in Theorem~\ref{thm:midlmc_complexity}. The slopes observed in Figures~\ref{fig:runtime} and~\ref{fig:work_est} match the predictions of Theorem~\ref{thm:midlmc_complexity}, lending credence to Assumptions~\ref{ass:midlmc_bias} and~\ref{ass:midlmc_var} for this problem. Figures~\ref{fig:runtime} and \ref{fig:work_est} support the superiority of the proposed multi-index DLMC estimator over the multilevel DLMC estimator, with the incorporated IS scheme enabling both estimators to estimate rare-event quantities feasibly. The multilevel DLMC failed to produce estimates of the quantity of interest up to $\tol_{\mathrm{r}} = 1\%$, whereas the multi-index algorithm achieved estimates within a fixed computational budget. Moreover, the multi-index estimator demonstrates a complexity reduction of about one order of magnitude (up to logarithmic terms) compared with the multilevel estimator, enabling the computation of rare-event quantities in MV-SDEs that were previously infeasible. The trade-off is the upfront computation of the IS control and the need for careful coupling of mixed-differences, but the results demonstrate that this overhead is justified for rare-events.

\section{Conclusion}
\label{sec:conclusion}

This work presents theoretical and numerical evidence under verifiable assumptions that the multi-index DLMC estimator outperforms the multilevel DLMC estimator~\citep{my_mldlmcis} for rare-event quantities in MV-SDEs. These quantities are expressed as expectations of sufficiently smooth observables of solutions to stochastic particle systems in the mean-field limit. This work applies the same  IS scheme in the work by~\citep{my_dlmcis} at every multi-index in the proposed multi-index DLMC estimator. The numerical results confirmed a significant variance reduction on all multi-indices. For the considered example, the multi-index DLMC estimator achieved a complexity of $\order{\tol_{\mathrm{r}}^{-2} (\log \tol_{\mathrm{r}}^{-1})^2}$, an order of improvement over the multilevel DLMC estimator, where $\tol_{\mathrm{r}}$ represents the relative error tolerance. Not only did the asymptotic exponent improve, but the constant prefactor is much smaller than that in a comparable MIMC method for a nonrare observable~\citep{mlmc_mvsde}. This improvement enables the computation of rare-event quantities with much less computational cost, enabling robust uncertainty quantification analysis for MV-SDEs, which was previously impractical.

Although this work presents numerical results for the Kuramoto model, the proposed method is quite general and readily extends to other MV-SDEs, such as animal flocking models~\citep{Cucker:2007aa} and social opinion dynamics~\citep{Degond:2016aa}. The critical requirements are sufficient regularity and a tractable numerical solution to the control PDE~\eqref{eqn:dmvsde_hjb_form3}. Similar computational gains are expected in these settings, although verifying Assumptions~\ref{ass:midlmc_bias} and~\ref{ass:midlmc_var} for each case is critical. Future work will focus on rigorously analyzing these mixed regularity assumptions to derive sufficient conditions on $b,\sigma,\kappa_1,\kappa_2$, and $G$. Extending the method to nonsmooth ($C^0$ or lower) rare-event observables  appears viable (as evidenced by the partial results with mollifiers; see Remark~\ref{rem:probabilities}), although an adaptive mollification or other variance control strategy may be necessary for discontinuous observables. Although the complexity bound~\eqref{eqn:midlmc_complexity} does not depend on $d$, solving the control PDE~\eqref{eqn:dmvsde_hjb_form3} in very high dimensions ($d > 10$) is impractical. In such cases, model reduction techniques or stochastic optimization could be explored to obtain an effective IS control (see Remark~\ref{rem:high_dim}). Finally, optimizing the multi-index DLMC parameters (e.g., by determining the optimal $\theta$~\citep{optimal_mlmc} or integrating a continuation-type multi-index algorithm~\citep{cont_mlmc}) could yield additional computational savings.

\clearpage
\appendix

\section{Double loop Monte Carlo algorithm with importance sampling}
\label{app:dlmc_outline}

\begin{algorithm}
    \caption{Double loop Monte Carlo algorithm with importance sampling for MV-SDEs}
	\label{alg:dlmc_outline}
    \SetAlgoLined
    \textbf{Off-line: } \\
    Generate one realization of law $\{X_p^{\bar{P}|\bar{N}}\}_{p=1}^{\bar{P}} \sim \mu^{\bar{P}|\bar{N}}$ with the $\bar{P}$-particle system and $\bar{N}$ time steps with a large value of $\bar{P},\bar{N}$; \\
    Given $\mu^{\bar{P}|\bar{N}}$, solve the PDE~\eqref{eqn:dmvsde_hjb_form3} to obtain IS control $\zeta$; \\
    \textbf{Inputs: } $P,N,M_1,M_2,\zeta$; \\
    \For{$i=1,\ldots,M_1$}{
    Generate realization of random variables $\omega_{1:P}^{(i)}$; \\
    Generate a realization of law $\{X_p^{P|N}\}_{p=1}^P \sim \mu^{P|N} \left( \omega_{1:P}^{(i)} \right)$ with the $P$-particle system and $N$ time steps; \\
    \For{$j=1,\ldots,M_2$}{
    Generate realization of random variables ${\tilde{\omega}^{(i,j)}}$; \\
    Given $\mu^{P|N} \left( \omega_{1:P}^{(i)} \right)$ and $\zeta$, generate a sample path of~\eqref{eqn:dmvsde_sde_is} with $N$ time steps;\\
    Compute $G \left(\Bar{X}_{\zeta}^{P|N}(T)\right) \left( \omega_{1:P}^{(i)} \times {\tilde{\omega}^{(i,j)}} \right)$; \\
    Compute $\mathbb{L}^{P|N} \left( \omega_{1:P}^{(i)} \times {\tilde{\omega}^{(i,j)}} \right)$ using~\eqref{eqn:dlmc_llhood_factor};
    }
    Approximate $\E{G \left(\Bar{X}_{\zeta}^{P|N}(T)\right) \mathbb{L}^{P|N} \mid \mu^{P|N} \left( \omega_{1:P}^{(i)} \right) }$ by $\frac{1}{M_2} \sum_{j=1}^{M_2} G \left(\Bar{X}_{\zeta}^{P|N}(T)\right) \mathbb{L}^{P|N} \left( \omega_{1:P}^{(i)} \times {\tilde{\omega}^{(i,j)}} \right)$; \\
    }
    Approximate $\E{G \left( \bar{X}^{P|N}(T) \right)}$ by $\frac{1}{M_1} \sum_{i=1}^{M_1} \frac{1}{M_2} \sum_{j=1}^{M_2}  G \left(\Bar{X}_{\zeta}^{P|N}(T)\right) \mathbb{L}^{P|N} \left( \omega_{1:P}^{(i)} \times {\tilde{\omega}^{(i,j)}} \right)$ ;
\end{algorithm}

\section{Proof of Theorem~\ref{thm:midlmc_complexity}}
\label{app:complexity}

For a given index set $\mathcal{I}$, the optimal $M_{1,\boldsymbol\alpha} \in \mathbb{R}_+$ and $M_{2,\boldsymbol\alpha} \in \mathbb{R}_+$ are determined for all $\boldsymbol\alpha \in \mathcal{I}$, minimizing the total computational cost~\eqref{eqn:midlmc_work} subject to the statistical error constraint~\eqref{eqn:midlmc_var_constraint}:

\begin{empheq}[left=\empheqlbrace, right = \cdot]{equation}
\label{eqn:midlmc_optimization}
		\begin{alignedat}{2}
			&\min_{\{M_{1,\boldsymbol\alpha},M_{2,\boldsymbol\alpha}\}_{\boldsymbol\alpha \in \mathcal{I}}} \sum_{\boldsymbol\alpha \in \mathcal{I}} M_{1,\boldsymbol\alpha} N_{\alpha_2} P_{\alpha_1}^{2} + M_{1,\boldsymbol\alpha} M_{2,\boldsymbol\alpha} N_{\alpha_2} P_{\alpha_1} \\ 
			&\text{s.t. } \left(\sum_{\boldsymbol\alpha \in \mathcal{I}} \frac{V_{1,\boldsymbol\alpha}}{M_{1,\boldsymbol\alpha}} + \frac{V_{2,\boldsymbol\alpha}}{M_{1,\boldsymbol\alpha} M_{2,\boldsymbol\alpha}} \right) \approx \left( \frac{\theta \tol_{\mathrm{r}} \E{G(X(T))}}{C_\nu} \right)^2
		\end{alignedat}
\end{empheq}
The solution to \eqref{eqn:midlmc_optimization} in $\mathbb{R}_+$ using the Lagrange multiplier method for constrained optimization yields the following:

\begin{align}
	&\mathcal{M}_{1,\boldsymbol\alpha} = \left( \frac{C_\nu}{\theta \tol_{\mathrm{r}} \E{G(X(T))}} \right)^2 \sqrt{\frac{V_{1,\boldsymbol\alpha}}{N_{\alpha_2} P_{\alpha_1}^{2}}} \sum_{\boldsymbol\beta \in \mathcal{I}} \left( \sqrt{V_{1,\boldsymbol\beta} N_{\beta_2} P_{\beta_1}^{2}} + \sqrt{V_{2,\boldsymbol\beta} N_{\beta_2} P_{\beta_1}} \right), \nonumber \\
	\label{eqn:midlmc_optimal_soln}
	&\tilde{\mathcal{M}}_{\boldsymbol\alpha} = \mathcal{M}_{1,\boldsymbol\alpha} \mathcal{M}_{2,\boldsymbol\alpha} \\
	&\qquad = \left( \frac{C_\nu}{\theta \tol_{\mathrm{r}} \E{G(X(T))}} \right)^2 \sqrt{\frac{V_{2,\boldsymbol\alpha}}{N_{\alpha_2} P_{\alpha_1}}} \sum_{\boldsymbol\beta \in \mathcal{I}} \left( \sqrt{V_{1,\boldsymbol\beta} N_{\beta_2} P_{\beta_1}^{2}} + \sqrt{V_{2,\boldsymbol\beta} N_{\beta_2} P_{\beta_1}} \right) \cdot \nonumber
\end{align}
In practice, natural numbers are employed for $\{\mathcal{M}_{1,\boldsymbol\alpha},\mathcal{M}_{2,\boldsymbol\alpha}\}_{\boldsymbol\alpha \in \mathcal{I}}$. Therefore, to guarantee at least one $\mathcal{M}_{1,\boldsymbol\alpha}$ and $\mathcal{M}_{2,\boldsymbol\alpha}$ for each $\boldsymbol\alpha$, the following quasi-optimal solution to~\eqref{eqn:midlmc_optimization} is employed:

\begin{equation}
	\label{eqn:midlmc_actual_soln}
		M_{1,\boldsymbol\alpha} = \lceil \mathcal{M}_{1,\boldsymbol\alpha} \rceil, \quad M_{2,\boldsymbol\alpha} = \left\lceil \frac{\tilde{\mathcal{M}}_{\boldsymbol\alpha}}{\lceil \mathcal{M}_{1,\boldsymbol\alpha} \rceil} \right\rceil \cdot
\end{equation}
Using \eqref{eqn:midlmc_actual_soln}, the computational cost of the estimator is bound as follows:

\begin{align}
		\mathcal{W}[\mathcal{A}_{\mathrm{MIMC}}(\mathcal{I})] &\lesssim \sum_{\boldsymbol\alpha \in \mathcal{I}} \left( (\mathcal{M}_{1,\boldsymbol\alpha}+1) N_{\alpha_2} P_{\alpha_1}^{2} + (\mathcal{M}_{1,\boldsymbol\alpha}+1) \left( \frac{\tilde{\mathcal{M}}_{\boldsymbol\alpha}}{\lceil \mathcal{M}_{1,\boldsymbol\alpha} \rceil}+1 \right) N_{\alpha_2} P_{\alpha_1} \right) \nonumber \\
		&\leq \underbrace{\sum_{\boldsymbol\alpha \in \mathcal{I}} \left( \mathcal{M}_{1,\boldsymbol\alpha} N_{\alpha_2} P_{\alpha_1}^{2} + \tilde{\mathcal{M}}_{\boldsymbol\alpha} N_{\alpha_2} P_{\alpha_1} \right)}_{=W_1(\mathcal{I})} \nonumber \\
		&\qquad + \underbrace{ \sum_{\boldsymbol\alpha \in \mathcal{I}} \left( P_{\alpha_1}^{2} N_{\alpha_2} + P_{\alpha_1} N_{\alpha_2} \right) }_{=W_2(\mathcal{I}),\text{ cost of one sample per multi-index}} \nonumber \\
		\label{eqn:midlmc_work_bound1}
		&\qquad + \underbrace{ \sum_{\boldsymbol\alpha \in \mathcal{I}} \mathcal{M}_{1,\boldsymbol\alpha} N_{\alpha_2} P_{\alpha_1}}_{=W_3(\mathcal{I})} + \underbrace{ \sum_{\boldsymbol\alpha \in \mathcal{I}} \frac{\tilde{\mathcal{M}}_{\boldsymbol\alpha}}{\lceil \mathcal{M}_{1,\boldsymbol\alpha} \rceil} N_{\alpha_2} P_{\alpha_1}}_{=W_4(\mathcal{I})} \cdot
\end{align}
Applying $P_{\alpha_1} \geq 1$ reveals that

\begin{align*}
	W_3(\mathcal{I}) &= \sum_{\boldsymbol\alpha \in \mathcal{I}} \mathcal{M}_{1,\boldsymbol\alpha} N_{\alpha_2} P_{\alpha_1} \leq \sum_{\boldsymbol\alpha \in \mathcal{I}} \mathcal{M}_{1,\boldsymbol\alpha} N_{\alpha_2} P_{\alpha_1}^{2} \leq W_1(\mathcal{I}), \\
	W_4(\mathcal{I}) &= \sum_{\boldsymbol\alpha \in \mathcal{I}} \frac{\tilde{\mathcal{M}}_{\boldsymbol\alpha}}{\lceil \mathcal{M}_{1,\boldsymbol\alpha} \rceil} N_{\alpha_2} P_{\alpha_1} \leq \sum_{\boldsymbol\alpha \in \mathcal{I}} \frac{\tilde{\mathcal{M}}_{\boldsymbol\alpha}}{\max(1,\mathcal{M}_{1,\boldsymbol\alpha})} N_{\alpha_2} P_{\alpha_1}\\
	&\qquad \leq \sum_{\boldsymbol\alpha \in \mathcal{I}} \tilde{\mathcal{M}}_{\boldsymbol\alpha} N_{\alpha_2} P_{\alpha_1} \leq W_1(\mathcal{I}) \cdot
\end{align*} 
Hence, \eqref{eqn:midlmc_work_bound1} is rewritten as follows:

\begin{equation}
	\label{eqn:midlmc_work_bound2}
	\mathcal{W}[\mathcal{A}_{\mathrm{MIMC}}(\mathcal{I})] \lesssim W_1(\mathcal{I}) + W_2(\mathcal{I}) \cdot
\end{equation}
Substituting \eqref{eqn:midlmc_optimal_soln} in \eqref{eqn:midlmc_work_bound2} yields

\begin{align}
	\label{eqn:midlmc_work_bound3}
	\mathcal{W}[\mathcal{A}_{\mathrm{MIMC}}(\mathcal{I})] &\lesssim \left( \frac{C_\nu}{\theta \tol_{\mathrm{r}} \E{G(X(T))}} \right)^2 \left( \sum_{\boldsymbol\alpha \in \mathcal{I}} \sqrt{V_{1,\boldsymbol\alpha} N_{\alpha_2} P_{\alpha_1}^{2}} + \sqrt{V_{2,\boldsymbol\alpha} N_{\alpha_2} P_{\alpha_1}} \right)^2 \\
	&\qquad + \sum_{\boldsymbol\alpha \in \mathcal{I}} \left( P_{\alpha_1}^{2} N_{\alpha_2} + P_{\alpha_1} N_{\alpha_2} \right) \cdot \nonumber 
\end{align}
Under Assumptions \ref{ass:midlmc_bias} and \ref{ass:midlmc_var}, the total computational cost \eqref{eqn:midlmc_work_bound3} to reach a relative error tolerance $\tol_\mathrm{r}$ is estimated as follows:

\begin{align}
		\label{eqn:midlmc_work_bound4}
		\mathcal{W}[\mathcal{A}_{\mathrm{MIMC}}(\mathcal{I})] &\lesssim \left( \frac{C_\nu}{\theta \tol_{\mathrm{r}} \E{G(X(T))}} \right)^2  \left( \sum_{\boldsymbol\alpha \in \mathcal{I}} \left( 2^{\frac{\alpha_1(2-w_1) + \alpha_2(1-w_2)}{2}} + 2^{\frac{\alpha_1(1-s_1) + \alpha_2(1-s_2)}{2}} \right) \right)^2 \\
		&\quad+ \sum_{\boldsymbol\alpha \in \mathcal{I}} 2^{2\alpha_1+\alpha_2} + 2^{\alpha_1+\alpha_2} \nonumber\\
		&\leq \left( \frac{2 C_\nu}{\theta \tol_{\mathrm{r}} \E{G(X(T))}} \right)^2  \left( \sum_{\boldsymbol\alpha \in \mathcal{I}} \left( 2^{\frac{\alpha_1(1-\min (w_1-1,s_1)) + \alpha_2(1-\min(w_2,s_2))}{2}} \right) \right)^2 \nonumber \\
		&\quad + 2 \sum_{\boldsymbol\alpha \in \mathcal{I}} 2^{2\alpha_1+\alpha_2} \cdot \nonumber
\end{align}
The computational cost of the estimator is bounded as follows by defining $\bar{s}_1=\min (w_1-1,s_1)$ and $\bar{s}_2=\min(w_2,s_2)$:

\begin{align}
	\label{eqn:midlmc_work_bound5}
	\mathcal{W}[\mathcal{A}_{\mathrm{MIMC}}(\mathcal{I})] &\lesssim \left( \frac{2 C_\nu}{\theta \tol_{\mathrm{r}} \E{G(X(T))}} \right)^2  \left( \underbrace{\sum_{\boldsymbol\alpha \in \mathcal{I}} \left( 2^{\frac{\alpha_1(1-\bar{s}_1) + \alpha_2(1-\bar{s}_2)}{2}} \right)}_{=\tilde{W}_1(\mathcal{I})} \right)^2 \\
	&\quad + 2 \underbrace{\sum_{\boldsymbol\alpha \in \mathcal{I}} 2^{2\alpha_1+\alpha_2}}_{=\tilde{W}_2(\mathcal{I})}, \nonumber
\end{align}
where $\tilde{W}_2(\mathcal{I})$ represents the total computational cost required to generate one sample for each multi-index $\boldsymbol\alpha \in \mathcal{I}$, which is the minimum cost for a DLMC estimator. It does not dominate the first term in the bound \eqref{eqn:midlmc_work_bound4}. 

Next, this work justifies the choice for the set of multi-indices $\mathcal{I}$ that minimizes $\mathcal{W}[\mathcal{A}_{\mathrm{MIMC}}(\mathcal{I})]$. First, the following bound is obtained for the relative bias of the multi-index DLMC estimator:

\begin{equation}
\label{eqn:midlmc_bias_constraint_v2}
	\epsilon_b (\mathcal{I}) = \frac{1}{\abs{\E{G(X(T))}}} \abs{\sum_{\boldsymbol\alpha \notin \mathcal{I}} \E{{\hat{\Delta} G_{\boldsymbol\alpha}}}} \leq \frac{1}{\abs{\E{G(X(T))}}} \sum_{\boldsymbol\alpha \notin \mathcal{I}} \abs{\E{{\hat{\Delta} G_{\boldsymbol\alpha}}}} \cdot
\end{equation}
Then, $\mathcal{I}$ is set to satisfy the following relative bias constraint:

\begin{equation}
\label{eqn:midlmc_bias_constraint_v3}
	\frac{1}{\abs{\E{G(X(T))}}} \sum_{\boldsymbol\alpha \notin \mathcal{I}} \abs{\E{{\hat{\Delta} G_{\boldsymbol\alpha}}}} \leq (1-\theta) \tol_{\mathrm{r}} \cdot
\end{equation}
Using Assumption~\ref{ass:midlmc_bias}, \eqref{eqn:midlmc_bias_constraint_v3} is written as follows:

\begin{equation}
	\label{eqn:midlmc_bias_bound}
	\tilde{B}(\mathcal{I}) := \sum_{\boldsymbol\alpha \notin \mathcal{I}} 2^{-\alpha_1 b_1 - \alpha_2 b_2} \leq \frac{(1-\theta) \tol_{\mathrm{r}} \abs{\E{G(X(T))}}}{Q_B} \cdot
\end{equation}
The solution to the following optimization problem to obtain the quasi-optimal set $\mathcal{I}$:

\begin{equation}
	\label{eqn:midlmc_set_optim_v2}
	\min_{\mathcal{I} \subset \mathbb{N}^2} \tilde{W}_1(\mathcal{I}) \quad \text{such that} \quad \tilde{B}(\mathcal{I}) \leq \frac{(1-\theta) \tol_{\mathrm{r}} \E{G(X(T))}}{Q_B} \cdot
\end{equation}
The method in the work by~\citep{mimc2016} demonstrates that~\eqref{eqn:midlmc_optimal_indexset} is the optimal index set in the sense of~\eqref{eqn:midlmc_set_optim_v2}. However, the conditions $2 b_1 \geq \bar{s}_1 - 1$ and $2 b_2 \geq \bar{s}_2 - 1$ must be satisfied to allow admissible index sets. Next, using the optimal index set~\eqref{eqn:midlmc_optimal_indexset} yields the following:

\begin{equation}
    \label{eqn:midlmc_work_bound_standard}
	\tilde{W}_1(\mathcal{I}(L)) = \sum_{\boldsymbol\alpha \in \mathbb{N}^2 : \boldsymbol\delta \cdot \boldsymbol\alpha \leq L} \exp (\mathbf{g} \cdot \boldsymbol\alpha),
\end{equation}
where $\mathbf{g} = \log 2 \left[ \frac{1-\bar{s}_1}{2}, \frac{1-\bar{s}_2}{2} \right]$ and $\boldsymbol\delta = \frac{\log 2}{C_\delta} \left[ \frac{1 - \bar{s}_1}{2} + b_1, \frac{1 - \bar{s}_2}{2} + b_2 \right]$.\\
In addition, $C_\delta = \log 2 \left(\frac{1 - \bar{s}_1}{2} + b_1 + \frac{1 - \bar{s}_2}{2} + b_2 \right)$ represents the normalizing constant, ensuring that $0 \leq \delta_1, \delta_2 \leq 1$ and $\delta_1 + \delta_2 = 1$. {For the optimal index set $\mathcal{I}(L)$ defined by~\eqref{eqn:midlmc_optimal_indexset}, computational cost bound~\eqref{eqn:midlmc_work_bound_standard} and the bias relation~\eqref{eqn:midlmc_bias_bound}, the asymptotic MIMC complexity estimates of~\citep{mimc2016} directly imply Theorem~\ref{thm:midlmc_complexity}, thus concluding the proof.} This computation is omitted here for ease of presentation.

\section{Algorithm~\ref{alg:adaptive_midlmcis} implementation details}
\label{app:algorithms}
Algorithm~\ref{alg:midlmcis_md} outlines the implementation of the IS scheme presented in Section~\ref{sec:midlmc} to estimate $\E{{\hat{\Delta} G_{\boldsymbol\alpha}}}$ for each $\boldsymbol\alpha \in \mathcal{I}$, which is necessary for the multi-index DLMC estimator~\eqref{eqn:midlmc_estimator}.

\begin{algorithm}
    \caption{Importance sampling scheme to estimate $\E{{\hat{\Delta} G_{\boldsymbol\alpha}}}$}
    \label{alg:midlmcis_md}
    \SetAlgoLined
    \textbf{Inputs: } $\boldsymbol\alpha,M_1,M_2,\zeta(\cdot,\cdot)$; \\
    \For{$m_1=1,\ldots,M_1$}{
    	Generate a realization of random variables $\omega_{1:P_{\alpha_1}}^{(\boldsymbol\alpha,m_1)}$; \\
    	Using $\omega_{1:P_{\alpha_1}}^{(\boldsymbol\alpha,m_1)}$, generate a realization of empirical law $\mu$ for all required discretizations, as per \eqref{eqn:antithetic_sampler} by generating sample paths of the particle system~\eqref{eqn:strong_approx_mvsde}; \\
        \For{$m_2=1,\ldots,M_2$}{
        	Generate a realization of random variables ${\tilde{\omega}^{(\boldsymbol\alpha,m_1,m_2)}}$; \\
            Using ${\tilde{\omega}^{(\boldsymbol\alpha,m_1,m_2)}}$, generate a sample path of \eqref{eqn:dmvsde_sde_is} for all required discretizations as per \eqref{eqn:antithetic_sampler} with control $\zeta$;\\
            Compute ${\hat{\Delta} G_{\boldsymbol\alpha}} (\omega_{1:P_{\alpha_1}}^{(\boldsymbol\alpha,m_1)} \times {\tilde{\omega}^{(\boldsymbol\alpha,m_1,m_2)}})$ using \eqref{eqn:antithetic_sampler}; \\
        }
        }
    Approximate $\E{{\hat{\Delta} G_{\boldsymbol\alpha}}}$ by $\frac{1}{M_1} \sum_{m_1=1}^{M_1} \frac{1}{M_2} \sum_{m_2=1}^{M_2} {\hat{\Delta} G_{\boldsymbol\alpha}} (\omega_{1:P_{\alpha_1}}^{(\boldsymbol\alpha,m_1)} \times {\tilde{\omega}^{(\boldsymbol\alpha,m_1,m_2)}})$; \\
\end{algorithm}

\subsection{Estimating the bias of index set $\mathcal{I}$}
\label{sec:adaptive_bias}

In Algorithm~\ref{alg:adaptive_midlmcis}, to determine the optimal index set $\mathcal{I}(L)$ using~\eqref{eqn:midlmc_optimal_indexset}, the rates $\{b_1,b_2\}$, $\{w_1,w_2\}$, and $\{s_1,s_2\}$ associated with Assumptions~\ref{ass:midlmc_bias} and \ref{ass:midlmc_var} must be estimated. Reliable pilot runs are necessary to obtain these estimates. However, such pilot runs are unnecessary for the adaptive DLMC algorithm~\citep{my_dlmcis} or multilevel DLMC algorithm~\citep{my_mldlmcis}. The optimal index set in~\eqref{eqn:midlmc_optimal_indexset} is defined using the parameter (threshold) $L$. Given that $\mathcal{I}(L) \subset \mathcal{I}(L+1)$, the following heuristic estimate~\citep{mimc2016} is applied for the absolute bias corresponding to index set $\mathcal{I}(L)$.

\begin{equation}
	\label{eqn:adaptive_midlmc_bias}
	\abs{\E{G - \mathcal{A}_{\mathrm{MIMC}}(\mathcal{I}(L))}} \approx \sum_{\boldsymbol\alpha \in \partial \mathcal{I}(L)} \abs{\E{{\Delta^\text{mix} G_{\boldsymbol\alpha}}}}, 
\end{equation}
where $\partial \mathcal{I}(L)$ denotes the outer boundary of the index set $\mathcal{I}(L)$, which is  defined as follows:

\begin{equation}
	\label{eqn:adaptive_midlmc_boundary}
	\partial \mathcal{I}(L) = \{\boldsymbol\alpha \in \mathcal{I}(L): \boldsymbol\alpha + (1,0) \notin \mathcal{I}(L) \quad \mathbf{or} \quad \boldsymbol\alpha + (0,1) \notin \mathcal{I}(L)\} \cdot
\end{equation}
Because $\partial \mathcal{I}(L) \subset \mathcal{I}(L)$, the computed nested averages are employed  using the optimal $\{M_{1,\boldsymbol\alpha},M_{2,\boldsymbol\alpha}\}_{\boldsymbol\alpha \in \partial \mathcal{I}(L)}$ in~\eqref{eqn:midlmc_estimator} to estimate each expectation in~\eqref{eqn:adaptive_midlmc_bias}. 

\subsection{Estimating $\{V_{1,\boldsymbol\alpha},V_{2,\boldsymbol\alpha}\}$}
\label{sec:adaptive_variance}

Computationally inexpensive, robust empirical estimates of $\{V_{1,\boldsymbol\alpha}, V_{2,\boldsymbol\alpha}\}$ for all $\boldsymbol\alpha \in \mathcal{I}(L)$ are necessary to compute the optimal number of samples to satisfy the variance constraint in~\eqref{eqn:midlmc_var_constraint} of the estimator using~\eqref{eqn:midlmc_actual_soln}. Therefore, the DLMC method in  Algorithm~\ref{alg:adaptive_variance} is employed with appropriately selected values of $\tilde{M}_1,\tilde{M}_2$.

\begin{algorithm}
    \caption{Estimating $\{V_{1,\boldsymbol\alpha},V_{2,\boldsymbol\alpha}\}$ for the adaptive multi-index double loop Monte Carlo method}
    \label{alg:adaptive_variance}
    \SetAlgoLined
    \textbf{Inputs: } $\boldsymbol\alpha,M_1,M_2,\zeta(\cdot,\cdot)$; \\
    \For{$m_1=1,\ldots,M_1$}{
    	Generate a realization of random variables $\omega_{1:P_{\alpha_1}}^{(\boldsymbol\alpha,m_1)}$; \\
    	Using $\omega_{1:P_{\alpha_1}}^{(\boldsymbol\alpha,m_1)}$, generate a realization of empirical law $\mu$ for all required discretizations, as per \eqref{eqn:antithetic_sampler} by generating sample paths of the particle system~\eqref{eqn:strong_approx_mvsde}; \\
        \For{$m_2=1,\ldots,M_2$}{
        	Generate a realization of random variables ${\tilde{\omega}^{(\boldsymbol\alpha,m_1,m_2)}}$; \\
            Using ${\tilde{\omega}^{(\boldsymbol\alpha,m_1,m_2)}}$, generate a sample path of \eqref{eqn:dmvsde_sde_is} for all required discretizations, as per \eqref{eqn:antithetic_sampler} with control $\zeta$;\\
            Compute ${\hat{\Delta} G_{\boldsymbol\alpha}} (\omega_{1:P_{\alpha_1}}^{(\boldsymbol\alpha,m_1)} \times {\tilde{\omega}^{(\boldsymbol\alpha,m_1,m_2)}})$ using \eqref{eqn:antithetic_sampler}; \\
        }
        Approximate $\E{{\hat{\Delta} G_{\boldsymbol\alpha}} \mid \omega_{1:P_{\alpha_1}}^{(\boldsymbol\alpha,m_1)}}$ by $\frac{1}{M_2} \sum_{m_2=1}^{M_2} {\hat{\Delta} G_{\boldsymbol\alpha}} (\omega_{1:P_{\alpha_1}}^{(\boldsymbol\alpha,m_1)} \times {\tilde{\omega}^{(\boldsymbol\alpha,m_1,m_2)}})$; \\
        Approximate $\Var{{\hat{\Delta} G_{\boldsymbol\alpha}} \mid \omega_{1:P_{\alpha_1}}^{(\boldsymbol\alpha,m_1)}}$ by the sample variance of $\left\{{\hat{\Delta} G_{\boldsymbol\alpha}} (\omega_{1:P_{\alpha_1}}^{(\boldsymbol\alpha,m_1)} \times {\tilde{\omega}^{(\boldsymbol\alpha,m_1,m_2)}})\right\}_{m_2=1}^{M_2}$; \\
        }
    Approximate $V_{1,\boldsymbol\alpha}$ by the sample variance of $\left\{\E{{\hat{\Delta} G_{\boldsymbol\alpha}} \mid \omega_{1:P_{\alpha_1}}^{(\boldsymbol\alpha,m_1)}}\right\}_{m_1=1}^{M_1}$; \\
    Approximate $V_{2,\boldsymbol\alpha}$ by $\frac{1}{M_1} \sum_{m_1=1}^{M_1} \Var{{\hat{\Delta} G_{\boldsymbol\alpha}} \mid \omega_{1:P_{\alpha_1}}^{(\boldsymbol\alpha,m_1)}} \cdot$ \\
\end{algorithm}

Estimating $\{V_{1,\boldsymbol\alpha},V_{2,\boldsymbol\alpha}\}$ for all $\boldsymbol\alpha \in \mathcal{I}(L)$ using Algorithm~\ref{alg:adaptive_variance} can become a computational burden when $\mathcal{I}(L)$ is even moderately large. To ease this computational overload, Assumption~\ref{ass:midlmc_var} is applied using an extrapolation approach to estimate $\{V_{1,\boldsymbol\alpha},V_{2,\boldsymbol\alpha}\}$ for deeper multi-indices.  Algorithm~\ref{alg:adaptive_variance} is applied solely to estimate $\{V_{1,\boldsymbol\alpha},V_{2,\boldsymbol\alpha}\}$ for the small full tensor index set $\{0,1,2\} \cross \{0,1,2\}$. Then, the extrapolation Algorithm~\ref{alg:midlmc_extrapolation} is employed using Assumption~\ref{ass:midlmc_var} to estimate $\{V_{1,\boldsymbol\alpha},V_{2,\boldsymbol\alpha}\}$ for the remaining multi-indices. To alleviate the computational burden further, this work only estimates $\{V_{1,\boldsymbol\alpha}, V_{2,\boldsymbol\alpha}\}$ for the newly added multi-indices in each iteration (i.e., for $\boldsymbol\alpha \in \mathcal{I}(L+1) - \mathcal{I}(L)$). As $\mathcal{I}(L) \subset \mathcal{I}(L+1)$ by construction, the estimates of $\{V_{1,\boldsymbol\alpha}, V_{2,\boldsymbol\alpha}\}$ can be reused for multi-indices from previous iterations.

\begin{algorithm}
    \caption{Sample variance extrapolation for the adaptive multi-index double loop Monte Carlo method}
\label{alg:midlmc_extrapolation}
    \SetAlgoLined
    \textbf{Inputs: } $\mathcal{I},\{\bar{M}_1,\bar{M}_2\}$,$(w_1,w_2)$,$(s_1,s_2)$; \\
    Estimate $\{V_{1,\boldsymbol\alpha},V_{2,\boldsymbol\alpha}\}$ for $\boldsymbol\alpha \in \{0,1,2\} \cross \{0,1,2\}$ using Algorithm~\ref{alg:adaptive_variance} with $\{\bar{M}_1,\bar{M}_2\}$; \\
    \For{$\boldsymbol\alpha \in \mathcal{I}$}
    {
    	\uIf{$\alpha_1=0,1$}{
    		$V_{1,\boldsymbol\alpha} = \max \left( \frac{V_{1,\boldsymbol\alpha-(0,1)}}{2^{w_2}}, \frac{V_{1,\boldsymbol\alpha-(0,2)}}{2^{2w_2}} \right)$;  $V_{2,\boldsymbol\alpha} = \max \left( \frac{V_{2,\boldsymbol\alpha-(0,1)}}{2^{s_2}}, \frac{V_{2,\boldsymbol\alpha-(0,2)}}{2^{2s_2}} \right)$; \\
    	}\uElseIf{$\alpha_2=0,1$}{
    	$V_{1,\boldsymbol\alpha} = \max \left( \frac{V_{1,\boldsymbol\alpha-(1,0)}}{2^{w_1}}, \frac{V_{1,\boldsymbol\alpha-(2,0)}}{2^{2w_1}} \right)$;  $V_{2,\boldsymbol\alpha} = \max \left( \frac{V_{2,\boldsymbol\alpha-(1,0)}}{2^{s_1}}, \frac{V_{2,\boldsymbol\alpha-(2,0)}}{2^{2s_1}} \right)$; \\
    	}\Else{
    	$V_{1,\boldsymbol\alpha} = \max \left( \frac{V_{1,\boldsymbol\alpha-(0,1)}}{2^{w_2}}, \frac{V_{1,\boldsymbol\alpha-(1,0)}}{2^{w_1}} \right)$; $V_{2,\boldsymbol\alpha} = \max \left( \frac{V_{2,\boldsymbol\alpha-(0,1)}}{2^{s_2}}, \frac{V_{2,\boldsymbol\alpha-(1,0)}}{2^{s_1}} \right)$; \\
    	} 
    }
\end{algorithm}

\subsection{Relative error control}

The adaptive algorithm requires a heuristic estimate of the quantity of interest $\E{G(X(T))}$ to meet the relative error constraints in~\eqref{eqn:midlmc_bias_constraint} and \eqref{eqn:midlmc_var_constraint} for a given relative error tolerance $\tol_{\mathrm{r}}$. In the algorithm, this estimate is updated at each iteration $L$. For $L=0$, Algorithm~\ref{alg:midlmcis_md} with appropriately selected $\bar{M}_1$ and $\bar{M}_2$ provides an initial estimate of $\bar{G}$ for $\E{G_{(0,0)}}$. In subsequent iterations, the multi-index estimator in~\eqref{eqn:midlmc_estimator} with optimal values of $\{M_{1,\boldsymbol\alpha},M_{2,\boldsymbol\alpha}\}_{\boldsymbol\alpha \in \mathcal{I}(L)}$ updates $\bar{G}$ and the absolute error tolerances. Combining these components, the adaptive multi-index DLMC Algorithm~\ref{alg:adaptive_midlmcis} is developed to evaluate rare-event observables associated with MV-SDEs. 

\section{Further numerical results}
\label{app:numerics_addendum}
\begin{figure}
    \centering
    \begin{subfigure}[b]{0.45\textwidth}
        \centering
        \includegraphics[width=\textwidth]{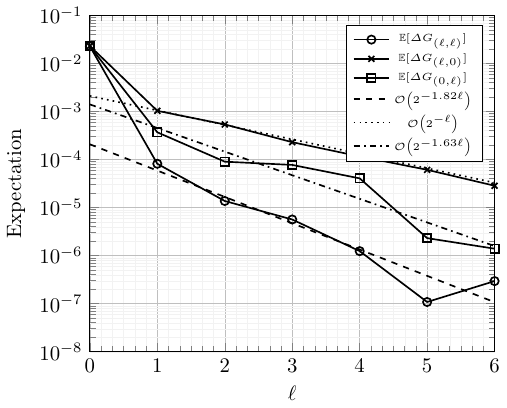}
        \caption{$C^0$-mollifier}
        \label{fig:c0_bias}
    \end{subfigure}
    \hfill
        \begin{subfigure}[b]{0.45\textwidth}
        \centering
        \includegraphics[width=\textwidth]{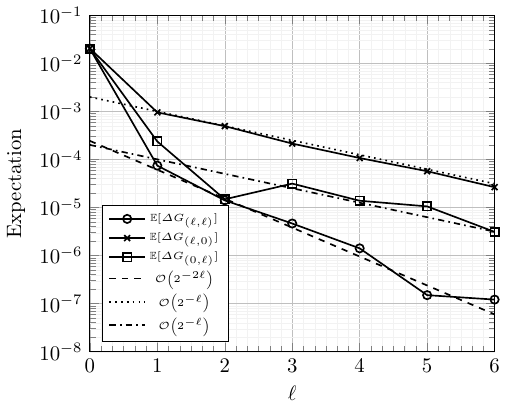}
        \caption{$C^1$-mollifier}
        \label{fig:c1_bias}
    \end{subfigure}
    \hfill
        \begin{subfigure}[b]{0.45\textwidth}
        \centering
        \includegraphics[width=\textwidth]{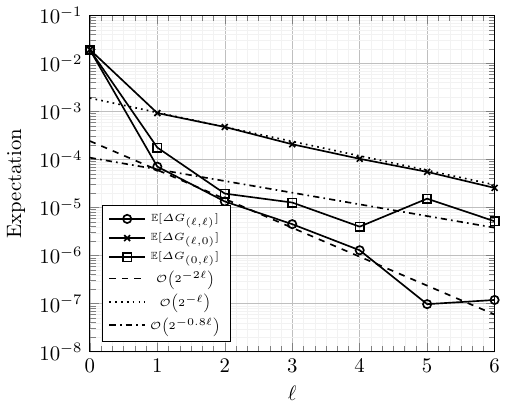}
        \caption{$C^2$-mollifier}
        \label{fig:c2_bias}
    \end{subfigure}
    \hfill
        \begin{subfigure}[b]{0.45\textwidth}
        \centering
        \includegraphics[width=\textwidth]{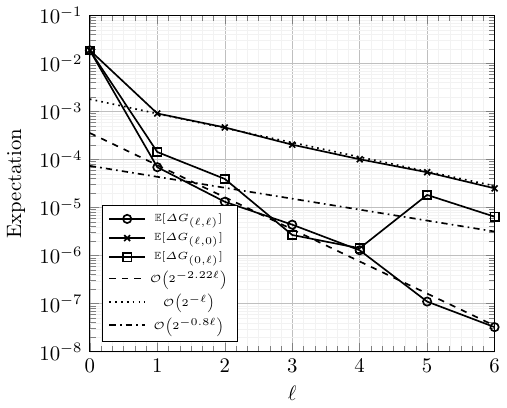}
        \caption{$C^3$-mollifier}
        \label{fig:c3_bias}
    \end{subfigure}
    \caption{Kuramoto example, numerical rate verification: Mixed-difference expectation convergence for mollifiers of varying regularity of the indicator function for $K=2$ and mollification parameter $\varepsilon=\frac{1}{2}$. Plots indicate that $C^1$-regularity suffices for multiplicative convergence of mixed-difference expectations, satisfying Assumption~\ref{ass:midlmc_bias}.}
    \label{fig:bias_addendum}
\end{figure}

\begin{figure}
    \centering
    \begin{subfigure}[b]{0.45\textwidth}
        \centering
        \includegraphics[width=\textwidth]{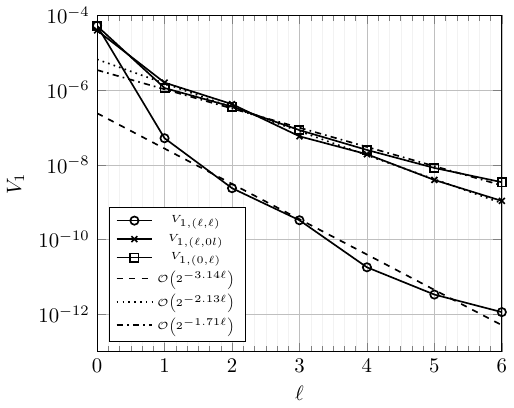}
        \caption{$C^0$-mollifier}
        \label{fig:c0_v1a}
    \end{subfigure}
    \hfill
        \begin{subfigure}[b]{0.45\textwidth}
        \centering
        \includegraphics[width=\textwidth]{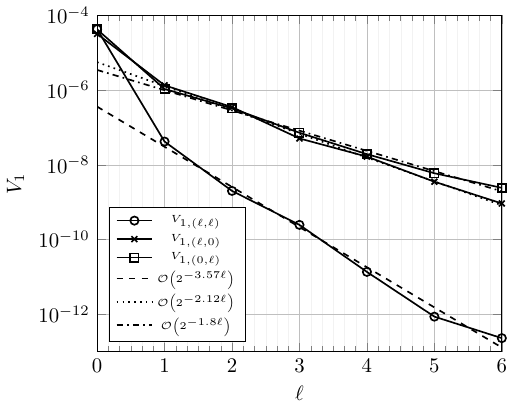}
        \caption{$C^1$-mollifier}
        \label{fig:c1_v1a}
    \end{subfigure}
    \hfill
        \begin{subfigure}[b]{0.45\textwidth}
        \centering
        \includegraphics[width=\textwidth]{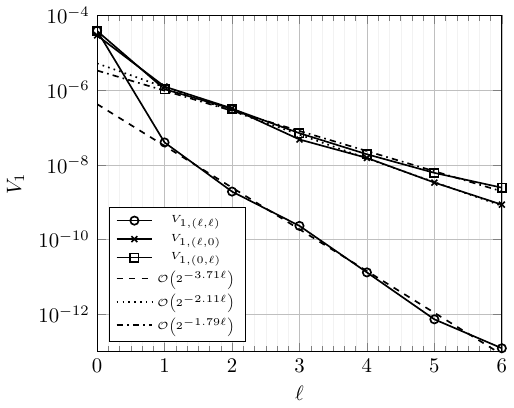}
        \caption{$C^2$-mollifier}
        \label{fig:c2_v1a}
    \end{subfigure}
    \hfill
        \begin{subfigure}[b]{0.45\textwidth}
        \centering
        \includegraphics[width=\textwidth]{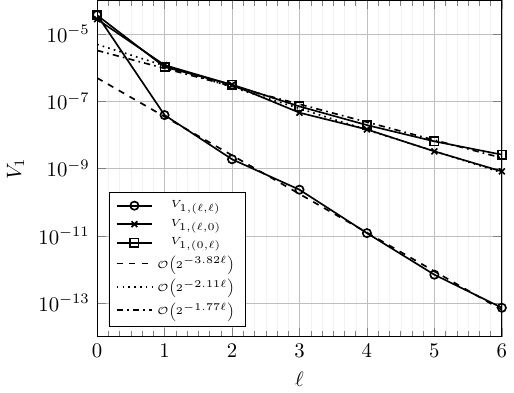}
        \caption{$C^3$-mollifier}
        \label{fig:c3_v1a}
    \end{subfigure}
    \caption{Kuramoto example, numerical rate verification: Mixed-difference $V_{1,\boldsymbol\alpha}$ convergence for mollifiers of varying regularity of the indicator function for $K=2$ and mollification parameter $\varepsilon=\frac{1}{2}$. Plots indicate that $C^2$-regularity suffices for multiplicative convergence of mixed-difference variances, satisfying Assumption~\ref{ass:midlmc_var}.}
    \label{fig:v1a_addendum}
\end{figure}

\begin{figure}
    \centering
    \begin{subfigure}[b]{0.45\textwidth}
        \centering
        \includegraphics[width=\textwidth]{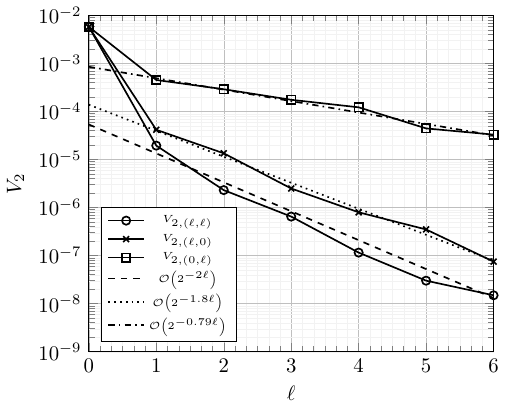}
        \caption{$C^0$-mollifier}
        \label{fig:c0_v2a}
    \end{subfigure}
    \hfill
        \begin{subfigure}[b]{0.45\textwidth}
        \centering
        \includegraphics[width=\textwidth]{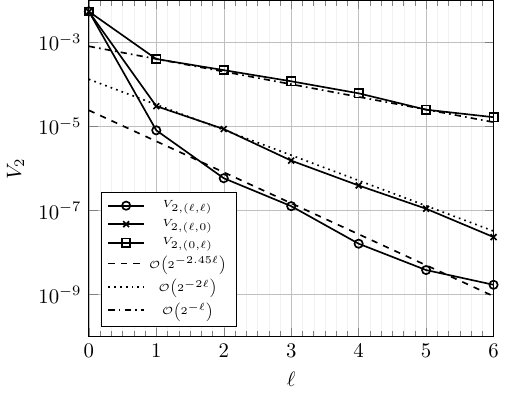}
        \caption{$C^1$-mollifier}
        \label{fig:c1_v2a}
    \end{subfigure}
    \hfill
        \begin{subfigure}[b]{0.45\textwidth}
        \centering
        \includegraphics[width=\textwidth]{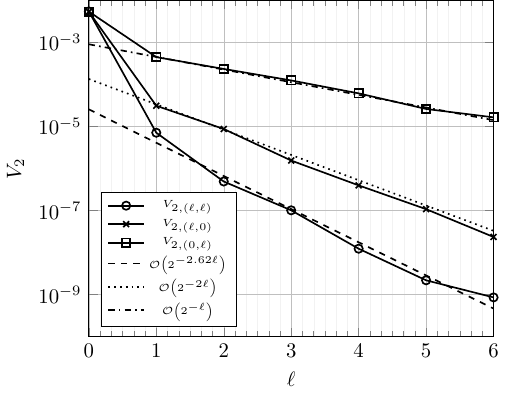}
        \caption{$C^2$-mollifier}
        \label{fig:c2_v2a}
    \end{subfigure}
    \hfill
        \begin{subfigure}[b]{0.45\textwidth}
        \centering
        \includegraphics[width=\textwidth]{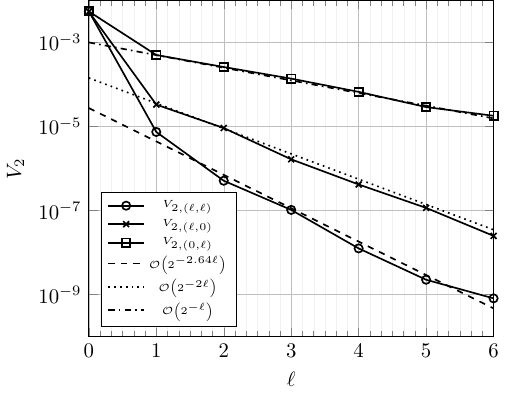}
        \caption{$C^3$-mollifier}
        \label{fig:c3_v2a}
    \end{subfigure}
    \caption{Kuramoto example, numerical rate verification: Mixed-difference $V_{2,\boldsymbol\alpha}$ convergence for mollifiers of varying regularity of the indicator function for $K=2$ and mollification parameter $\varepsilon=\frac{1}{2}$. Plots indicate that $C^2$-regularity suffices for multiplicative convergence of mixed-difference variances, satisfying Assumption~\ref{ass:midlmc_var}.}
    \label{fig:v2a_addendum}
\end{figure}

Figures~\ref{fig:bias_addendum},\ref{fig:v1a_addendum}, and~\ref{fig:v2a_addendum} corroborate the numerical rates depicted in Table~\ref{tab:regularity_rates}. The corresponding plots for the indicator function are omitted because they were erratic, which is expected, as a discontinuous observable does not satisfy the necessary regularity for Assumption~\ref{ass:midlmc_var}. This work also refrains from producing the corresponding plots for the Kuramoto model with a less regular triangular wave kernel, although these plots can be produced upon request.

\clearpage
{
\textbf{Declaration of generative AI and AI-assisted technologies in the manuscript preparation process}: During the preparation of this work, the authors used Anthropic Claude and OpenAI ChatGPT as editing tools. The authors reviewed and edited the output as needed and take full responsibility for the content of the published article. 
}
  \bibliographystyle{elsarticle-num} 
  \bibliography{references.bib}

@article{mimc2016,
	title={Multi-index {M}onte {C}arlo: when sparsity meets sampling},
  author={Haji-Ali, Abdul-Lateef and Nobile, Fabio and Tempone, Ra{\'u}l},
  journal={Numerische Mathematik},
  volume={132},
  pages={767--806},
  year={2016},
  publisher={Springer}}

@article{my_mldlmcis,
	author = {Ben Rached, Nadhir and Haji-Ali, Abdul-Lateef and Subbiah Pillai, Shyam Mohan and Tempone, Ra{\'u}l},
	journal = {Statistics and Computing},
	number = {1},
	pages = {1},
	title = {Multilevel importance sampling for rare events associated with the {M}c{K}ean--{V}lasov equation},
	volume = {35},
	year = {2024}}

@article{my_dlmcis,
	author = {Ben Rached, Nadhir and Haji-Ali, Abdul-Lateef and Subbiah Pillai, Shyam Mohan and Tempone, Ra{\'u}l},
	journal = {Statistics and Computing},
	number = {6},
	pages = {197},
	title = {Double-loop importance sampling for {M}c{K}ean--{V}lasov stochastic differential equation},
	volume = {34},
	year = {2024}}

@article{theory_mvsde,
	author = {Haji-Ali, Abdul-Lateef and Hoel, H{\aa}kon and Tempone, Ra{\'u}l},
	journal = {IMA Journal of Applied Mathematics},
	month = {9/10/2025},
	pages = {hxaf015},
	title = {Weak convergence analysis in the particle limit of the {McKean--Vlasov} equation using stochastic flows of particle systems},
	year = {2025}}

@article{chemical_app,
	title={The {K}uramoto model: A simple paradigm for synchronization phenomena},
  author={Acebr{\'o}n, Juan A and Bonilla, Luis L and Vicente, Conrad J P{\'e}rez and Ritort, F{\'e}lix and Spigler, Renato},
  journal={Reviews of modern physics},
  volume={77},
  number={1},
  pages={137},
  year={2005},
  publisher={APS}}

@article{neuroscience_app,
	title={Generalising the {K}uramoto model for the study of neuronal synchronisation in the brain},
  author={Cumin, David and Unsworth, CP},
  journal={Physica D: Nonlinear Phenomena},
  volume={226},
  number={2},
  pages={181--196},
  year={2007},
  publisher={Elsevier}}

@article{combustion_app,
	title={Arrays of coupled chemical oscillators},
  author={Forrester, Derek Michael},
  journal={Scientific reports},
  volume={5},
  number={1},
  pages={16994},
  year={2015},
  publisher={Nature Publishing Group UK London}}

@article{mlmc_mvsde,
	title={Multilevel and Multi-index {M}onte {C}arlo methods for the {M}c{K}ean--{V}lasov equation},
  author={Haji-Ali, Abdul-Lateef and Tempone, Ra{\'u}l},
  journal={Statistics and Computing},
  volume={28},
  pages={923--935},
  year={2018},
  publisher={Springer}}

@article{is_mvsde,
	title={Importance sampling for {M}c{K}ean-{V}lasov {SDE}s},
  author={dos Reis, Gon{\c{c}}alo and Smith, Greig and Tankov, Peter},
  journal={Applied Mathematics and Computation},
  volume={453},
  pages={128078},
  year={2023},
  publisher={Elsevier}}

@book{is_general_ref,
	title={Handbook of {M}onte {C}arlo methods},
  author={Kroese, Dirk P and Taimre, Thomas and Botev, Zdravko I},
  year={2013},
  publisher={John Wiley \& Sons}}

@article{mean_field_limit,
	title={Topics in propagation of chaos},
  author={Sznitman, Alain-Sol},
  journal={Lecture notes in mathematics},
  pages={165--251},
  year={1991},
  publisher={Springer Berlin Heidelberg}}

@book{sde_numerics,
	Author = {Kloeden, Peter E. and Platen, Eckhard},
	Title = {Numerical Solution of Stochastic Differential Equations},
	Publisher = {Springer, Berlin},
	doi = {https://doi.org/10.1007/978-3-662-12616-5},
	Year = {1992}}

@article{mckean_vlasov,
	title={A class of {M}arkov processes associated with nonlinear parabolic equations},
  author={McKean Jr, Henry P},
  journal={Proceedings of the National Academy of Sciences},
  volume={56},
  number={6},
  pages={1907--1911},
  year={1966},
  publisher={National Acad Sciences}}

@article{pedestrian_app,
	author = {Helbing, Dirk and Moln{\'a}r, P{\'e}ter},
	journal = {Physical Review E},
	month = {05},
	number = {5},
	pages = {4282--4286},
	title = {Social force model for pedestrian dynamics},
	volume = {51},
	year = {1995}}

@article{animal_app,
	author = {Erban, Radek and Haskovec, Jan},
	journal = {Kinetic and Related Models},
	number = {4},
	pages = {817--842},
	title = {From individual to collective behaviour of coupled velocity jump processes: A locust example},
	volume = {5},
	year = {2012}}

@article{biology_app,
	title={Particle-based multiscale modeling of calcium puff dynamics},
  author={Dobramysl, Ulrich and R{\"u}diger, Sten and Erban, Radek},
  journal={Multiscale Modeling \& Simulation},
  volume={14},
  number={3},
  pages={997--1016},
  year={2016},
  publisher={SIAM}}

@article{finance_app,
	 title={Stochastic evolution equations in portfolio credit modelling},
  author={Bush, Nick and Hambly, Ben M and Haworth, Helen and Jin, Lei and Reisinger, Christoph},
  journal={SIAM Journal on Financial Mathematics},
  volume={2},
  number={1},
  pages={627--664},
  year={2011},
  publisher={SIAM}}

@article{mvsde_existence,
	title={Existence and uniqueness theorems for solutions of {M}c{K}ean--{V}lasov stochastic equations},
  author={Mishura, Yuliya and Veretennikov, Alexander},
  journal={Theory of Probability and Mathematical Statistics},
  volume={103},
  pages={59--101},
  year={2020}}

@article{mvsde_pde,
	author = {Buckdahn, Rainer and Li, Juan and Peng, Shige and Rainer, Catherine},
	journal = {The Annals of Probability},
	month = {2025/11/28/},
	number = {2},
	pages = {824--878},
	title = {MEAN-FIELD STOCHASTIC DIFFERENTIAL EQUATIONS AND ASSOCIATED {PDE}S},
	volume = {45},
	year = {2017}}

@article{euler_mvsde,
  title={Strong convergence of {E}uler--{M}aruyama schemes for {M}c{K}ean--{V}lasov stochastic differential equations under local {L}ipschitz conditions of state variables},
  author={Li, Yun and Mao, Xuerong and Song, Qingshuo and Wu, Fuke and Yin, George},
  journal={IMA Journal of Numerical Analysis},
  volume={43},
  number={2},
  pages={1001--1035},
  year={2023},
  publisher={Oxford University Press}
}

@article{optimal_mlmc,
	title={Optimization of mesh hierarchies in multilevel {M}onte {C}arlo samplers},
  author={Haji-Ali, Abdul-Lateef and Nobile, Fabio and von Schwerin, Erik and Tempone, Ra{\'u}l},
  journal={Stochastics and Partial Differential Equations Analysis and Computations},
  volume={4},
  number={1},
  pages={76--112},
  year={2016},
  publisher={Springer}}

@article{cont_mlmc,
	title={A continuation multilevel {M}onte {C}arlo algorithm},
  author={Collier, Nathan and Haji-Ali, Abdul-Lateef and Nobile, Fabio and Von Schwerin, Erik and Tempone, Ra{\'u}l},
  journal={BIT Numerical Mathematics},
  volume={55},
  pages={399--432},
  year={2015},
  publisher={Springer}}

@article{mlmcis_kebaier,
	title={Coupling importance sampling and multilevel {M}onte {C}arlo using sample average approximation},
  author={Kebaier, Ahmed and Lelong, J{\'e}r{\^o}me},
  journal={Methodology and Computing in Applied Probability},
  volume={20},
  pages={611--641},
  year={2018},
  publisher={Springer}}

@article{mlmcis_alaya,
	title={Adaptive importance sampling for multilevel {M}onte {C}arlo {E}uler method},
  author={Ben Alaya, Mohamed and Hajji, Kaouther and Kebaier, Ahmed},
  journal={Stochastics},
  volume={95},
  number={2},
  pages={303--327},
  year={2023},
  publisher={Taylor \& Francis}}

@article{mlmcis_giles,
	title={Multilevel {M}onte {C}arlo method for ergodic {SDE}s without contractivity},
  author={Fang, Wei and Giles, Michael B},
  journal={Journal of Mathematical Analysis and Applications},
  volume={476},
  number={1},
  pages={149--176},
  year={2019},
  publisher={Elsevier}}

@article{mlmcis_srn,
	title={Importance sampling for a robust and efficient multilevel {M}onte {C}arlo estimator for stochastic reaction networks},
  author={Ben Hammouda, Chiheb and Ben Rached, Nadhir and Tempone, Ra{\'u}l},
  journal={Statistics and Computing},
  volume={30},
  pages={1665--1689},
  year={2020},
  publisher={Springer}}

@article{mvsde_weak_soln,
	author = {William R. P. Hammersley and David {\v S}i{\v s}ka and {\L}ukasz Szpruch},
	journal = {The Annals of Probability},
	month = {3},
	number = {2},
	pages = {527--555 },
	title = {Weak existence and uniqueness for McKean--Vlasov SDEs with common noise},
	volume = {49},
	year = {2021}}

@article{mvsde_strong_conv_1,
  title={A stochastic particle method for the {M}c{K}ean-{V}lasov and the {B}urgers equation},
  author={Bossy, Mireille and Talay, Denis},
  journal={Mathematics of computation},
  volume={66},
  number={217},
  pages={157--192},
  year={1997}
}

@article{mvsde_strong_conv_2,
  title={Convergence rate for the approximation of the limit law of weakly interacting particles: application to the {B}urgers equation},
  author={Bossy, Mireille and Talay, Denis},
  journal={The Annals of Applied Probability},
  volume={6},
  number={3},
  pages={818--861},
  year={1996},
  publisher={Institute of Mathematical Statistics}
}

@incollection{mvsde_strong_conv_3,
  Title={Asymptotic Behaviour of some interacting particle systems; {M}c{K}ean-{V}lasov and {B}oltzmann models},
  author={M\'el\'eard,S.},
  editor={Talay,D. and Tubaro,L.},
  booktitle={Probabilistic Models for Nonlinear Partial Differential Equations},
  publisher={Springer},
  year={1996},
  volume={1627},
  pages = {42--95}
}

@article{talay1990expansion,
	author = {Talay, Denis and Tubaro, Luciano},
	journal = {Stochastic Analysis and Applications},
	month = {01},
	number = {4},
	pages = {483--509},
	title = {Expansion of the global error for numerical schemes solving stochastic differential equations},
	volume = {8},
	year = {1990}}

@article{pflaum1999error,
  title={Error analysis of the combination technique},
  author={Pflaum, Christoph and Zhou, Aihui},
  journal={Numerische Mathematik},
  volume={84},
  number={2},
  pages={327--350},
  year={1999},
  publisher={Springer}
}

@proceedings{giles2023mlmc,
	author={Giles, Michael B.},
	address = {Cham},
	editor = {Hinrichs, Aicke and Kritzer, Peter and Pillichshammer, Friedrich},
	publisher = {Springer International Publishing},
	title = {MLMC Techniques for Discontinuous Functions},
	year = {2024}
}

@article{giles2015distribution,
author = {Giles, Michael B. and Nagapetyan, Tigran and Ritter, Klaus},
title = {Multilevel {M}onte {C}arlo Approximation of Distribution Functions and Densities},
journal = {SIAM/ASA Journal on Uncertainty Quantification},
volume = {3},
number = {1},
pages = {267-295},
year = {2015},
doi = {10.1137/140960086},
}

@article{szepessy2001adaptive,
author = {Szepessy, Anders and Tempone, Raúl and Zouraris, Georgios E.},
title = {Adaptive weak approximation of stochastic differential equations},
journal = {Communications on Pure and Applied Mathematics},
volume = {54},
number = {10},
pages = {1169-1214},
year = {2001}
}

@article{Hartmann:2016aa,
	author = {Hartmann, Carsten and Sch{\"u}tte, Christof and Zhang, Wei},
	journal = {Nonlinearity},
	number = {8},
	pages = {2298},
	title = {Model reduction algorithms for optimal control and importance sampling of diffusions},
	volume = {29},
	year = {2016}}

@article{Hartmann:2019aa,
	author = {Hartmann, Carsten and Kebiri, Omar and Neureither, Lara and Richter, Lorenz},
	journal = {Chaos: An Interdisciplinary Journal of Nonlinear Science},
	month = {2/20/2025},
	number = {6},
	pages = {063107},
	title = {Variational approach to rare event simulation using least-squares regression},
	volume = {29},
	year = {2019}}

@article{Khoromskij:2015aa,
	author = {Khoromskij, Boris N. },
	journal = {ESAIM: Proc.},
	month = {1},
	pages = {1--28},
	title = {Tensor numerical methods for multidimensional {PDE}S: theoretical analysis and initial applications},
	volume = {48},
	year = {2015}}

@article{Han:2018aa,
	author = {Han, Jiequn and Jentzen, Arnulf and E, Weinan},
	journal = {Proceedings of the National Academy of Sciences},
	month = {2025/09/24},
	number = {34},
	pages = {8505--8510},
	title = {Solving high-dimensional partial differential equations using deep learning},
	volume = {115},
	year = {2018}}

@article{Chaudru-de-Raynal:2020aa,
	author = {Chaudru de Raynal, P. E.},
	journal = {Stochastic Processes and their Applications},
	number = {1},
	pages = {79--107},
	title = {Strong well posedness of McKean--Vlasov stochastic differential equations with H{\"o}lder drift},
	volume = {130},
	year = {2020}}

@article{Chaudru-de-Raynal:2022aa,
	author = {Chaudru de Raynal, Paul-Eric and Frikha, Noufel},
	journal = {Journal de Math{\'e}matiques Pures et Appliqu{\'e}es},
	pages = {1--167},
	title = {Well-posedness for some non-linear SDEs and related PDE on the Wasserstein space},
	volume = {159},
	year = {2022}}

@article{Mischler:2015aa,
	author = {Mischler, St{\'e}phane and Mouhot, Cl{\'e}ment and Wennberg, Bernt},
	journal = {Probability Theory and Related Fields},
	number = {1},
	pages = {1--59},
	title = {A new approach to quantitative propagation of chaos for drift, diffusion and jump processes},
	volume = {161},
	year = {2015}}

@article{Bencheikh:2019aa,
	author = {Bencheikh, O. and Jourdain, B.},
	journal = {ESAIM: ProcS},
	pages = {219--235},
	title = {Bias behaviour and antithetic sampling in mean-field particle approximations of {SDE}s nonlinear in the sense of {McK}ean},
	volume = {65},
	year = {2019}}

@article{Szpruch:2021aa,
	author = {Szpruch, {\L}ukasz and Tse, Alvin},
	journal = {The Annals of Applied Probability},
	month = {2025/11/20/},
	number = {3},
	pages = {1100--1139},
	title = {ANTITHETIC MULTILEVEL SAMPLING METHOD FOR NONLINEAR FUNCTIONALS OF MEASURE},
	volume = {31},
	year = {2021}}

@article{Chassagneux:2022aa,
	author = {Chassagneux, Jean-Fran{\c c}ois and Szpruch, Lukasz and Tse, Alvin},
	journal = {The Annals of Applied Probability},
	number = {3},
	pages = {1929--1969},
	title = {Weak quantitative propagation of chaos via differential calculus on the space of measures},
	volume = {32},
	year = {2022}}

@article{Cucker:2007aa,
	author = {Cucker, Felipe and Smale, Steve},
	journal = {Japanese Journal of Mathematics},
	number = {1},
	pages = {197--227},
	title = {On the mathematics of emergence},
	volume = {2},
	year = {2007}}

@article{Degond:2016aa,
	author = {Degond, Pierre and Liu, Jian-Guo and Merino-Aceituno, Sara and Tardiveau, Thomas},
	journal = {Mathematical Models and Methods in Applied Sciences},
	month = {2025/11/28},
	number = {01},
	pages = {159--182},
	title = {Continuum dynamics of the intention field under weakly cohesive social interaction},
	volume = {27},
	year = {2016}}

@article{BenAmar:2026aa,
      title={Hierarchical Importance Sampling for Estimating Occupation Time for {SDE} Solutions}, 
      author={Ben Amar, Eya and Ben Rached, Nadhir and Tempone, Raúl},
      year={2025},
      journal = {arXiv preprint arXiv:2509.13950}, 
}



%
%
%
\end{document}